\documentclass[english,11pt]{article}

\usepackage{amsthm}
\usepackage[latin9]{inputenc}
\usepackage{babel}
\usepackage[hmargin=0.9in,vmargin=1.25in]{geometry}
\usepackage{subfig}
\usepackage{graphicx}
\usepackage[colorlinks=true,citecolor=blue,urlcolor=blue]{hyperref}

\usepackage{float}
\usepackage{url}
\usepackage{bm}

\usepackage{amsmath}
\usepackage{amssymb,amsfonts}
\usepackage{upgreek}
\usepackage[textsize=footnotesize]{todonotes}
\usepackage{stmaryrd}

\usepackage{mwe,tikz}\usepackage[percent]{overpic}

\newtheorem{remark}{Remark}
\numberwithin{remark}{subsection}

\usepackage[export]{adjustbox}



\graphicspath{{./figures/}}

\newcommand{\order}{{\mathcal O}}

\newcommand{\xihat}{\hat{\xi}}

\newcommand{\rhostar}{\overline{\rho}(\theta)}
\newcommand{\Kbar}{\overline{K}}
\newcommand{\sigmahat}{\hat{\sigma}}
\newcommand{\sigmabar}{\overline{\sigma}}
\def\avg[#1]{\left\langle #1 \right\rangle}
\def\intavg[#1]{\llbracket #1 \rrbracket}

\newcommand{\bfx}{{\bf x}}

\begin{document}
\title{Effective Rankine-Hugoniot conditions for shock waves in periodic media}
\author{
  David I. Ketcheson\thanks{Computer, Electrical, and Mathematical Sciences \& Engineering Division,
King Abdullah University of Science and Technology, 4700 KAUST, Thuwal
23955, Saudi Arabia. (david.ketcheson@kaust.edu.sa)} \and
  Manuel Quezada de Luna\thanks{Computer, Electrical, and Mathematical Sciences \& Engineering Division, 
King Abdullah University of Science and Technology, 4700 KAUST, Thuwal
23955, Saudi Arabia. (manuel.quezada@kaust.edu.sa.)}}
\maketitle

\begin{abstract}
  Solutions of first-order nonlinear hyperbolic conservation laws typically
  develop shocks in finite time even from smooth initial conditions.
  However, in heterogeneous media with rapid spatial variation, shock formation may be delayed or avoided.
  When shocks do form in such media, their speed of propagation depends
  on the material structure.  We investigate conditions for
  shock formation and propagation in heterogeneous media.
  We focus on the propagation of plane waves in two-dimensional
  media with a periodic structure that changes in only one direction.
  %periodic media with material variation in only one direction.
  %
  We propose an estimate for the speed of the shocks that is based on the Rankine-Hugoniot
  conditions applied to a leading-order homogenized (constant coefficient) system.
  We verify this estimate via numerical simulations using different nonlinear constitutive relations
  and layered and smoothly varying media with a periodic structure.
  %periodic media.
  %
  In addition, we discuss conditions and regimes under which shocks form in this type of media. 
\end{abstract}

%\begin{keywords}  shock wave; periodic medium; dispersion; homogenization
%\end{keywords}

% \begin{AMS} 35B27; 35L60; 35L67
%\end{AMS}

\section{Introduction}
Many important physical phenomena are modeled by first-order nonlinear
hyperbolic conservation laws; important examples include water waves and fluid
dynamics.  It is well known that solutions of such equations generically
develop singularities (shocks) in finite time and eventually dissipate all
energy.  This was shown rigorously for unbounded homogeneous domains with
compact initial data in \cite{liu1977decay, liu1977linear}.

  There are, however, mechanisms that can impede the formation of shocks.
  For instance, the effect of random topography on shallow water waves was
  studied in \cite{fouque2004shock,garnier2007effective}, where it is
  shown to lead to an effective viscosity.  As a result, waves over
  random topography propagate over longer distances before breaking.

For gas dynamics in a periodic domain, there seem to be solutions that remain
regular for long times; see \cite{shefter1999,qu2015long}.
Indeed, in \cite{shefter1999} it was conjectured, based on numerical
experiments, that in a periodic domain there exist nontrivial solutions in
which no singularity forms.  These solutions were referred to as NBAT
(non-breaking for all time).

  One may also consider a setting in which the medium properties vary periodically
  in space but the solution itself is not assumed periodic.  This is what is meant
  by the phrase {\em periodic medium} in the present work.
This situation arises naturally in applications that include
photonic \cite{halevi1999photonic} and phononic
crystals \cite{kushwaha1993acoustic, pennec2010two} and even coastal engineering
\cite{davies1984surface, tang2006omnidirectional}.
  The periodic variation in the medium leads to an effective
  dispersion \cite{santosa1991} for linear waves.
%
%If the problem is nonlinear, this dispersion can
%act as a regularization mechanism that prevents the formation of singularities.
%Therefore, NBAT solutions can also exist in periodic media.
%
Nonlinear waves in a periodic medium were studied in \cite{leveque2003}, where the authors 
derived via asymptotic expansions an effecitve constant coefficient system of PDEs,
whose leading order terms capture the macroscopic behavior
of the solution. The lowest-order terms match those of the original system, but with
spatially averaged coefficients,
while the higher-order terms are dispersive.

\begin{figure}[ht]
\begin{centering}
  \includegraphics[scale=0.25]{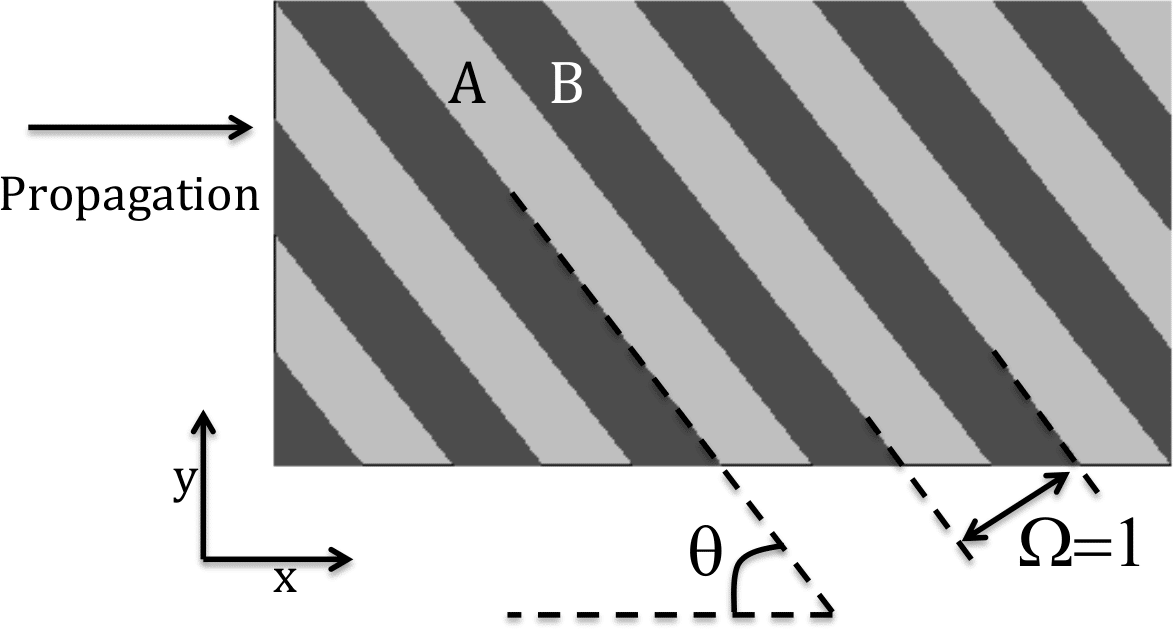}
\par\end{centering}
\caption{\small Wave propagation in a layered medium in two dimensions. 
  The domain shown is repeated periodically and extends infinitely in both directions.
  The orientation of the medium relative to the initial plane wave is given by the angle $\theta$.
  If $\theta=90^{\circ}$ we refer to the propagation as {\em transverse} to the material
  heterogeneity and if $\theta=0^{\circ}$ we refer to the propagation as {\em parallel}.
  \label{fig:layered_intro}}
\end{figure}

In \cite{shefter1999} and \cite{fouque2004shock,garnier2007effective,santosa1991,leveque2003},
regularization seems to be the result of reflection; in the first
case due to resonant interactions between characteristic fields and in the rest
due to reflections caused by the medium itself, which lead to dissipative and/or
dispersive effects.
In all these references, the authors focused on waves in one space dimension. 
In the present work we consider an unbounded two-dimensional domain that
varies periodically in one direction;
see Figure \ref{fig:layered_intro}.
This scenario has been studied already in \cite{quezada2014two} for linear waves and in
\cite{ketcheson2015diffractons} for nonlinear waves.
In the last reference, a different sort of regularization was observed that
does not seem to be caused by reflection, showing that apparently NBAT
solutions can arise even when reflective effects are essentially absent.

These mechanisms that prevent or delay the formation of shocks in
\cite{fouque2004shock,garnier2007effective,leveque2003,ketcheson2015diffractons}
depend on the initial data and the properties of the medium (i.e., the strengh of the
changes in the coefficients, topography, etc.). It is, therefore, possible
for shocks to still appear if the induced regularization is not strong enough.
Two natural questions arise:
\begin{itemize}
\item[i)] What properties of the initial data and the medium determine whether
  solutions exhibit shocks or are NBAT?
\item[ii)] When shocks form, what is their speed of propagation?
\end{itemize}
For one-dimensional layered media, this was studied in \cite{Ketcheson_LeVeque_2011},
where it was conjectured (based on numerical experiments) that NBAT solutions
also exist. The maximum amplitude of initial data leading to NBAT solutions there
was related to the properties of the medium. It was also demonstrated numerically
that initially discontinuous data could become regular (or almost regular) in such
media if the initial discontinuities did not satisfy a certain effective Lax-entropy condition.
  For shocks that do satisfy this condition, the authors proposed (and verified numerically)
  a simple expression for the speed of propagation.

Our goal is to extend the results in \cite{Ketcheson_LeVeque_2011} to the more
general two-dimensional situation depicted in Figure \ref{fig:layered_intro}.
To answer question (ii), we assume that the regularization is not strong enough, which leads
to the formation of shocks. We propose and test numerically a simple expression for the
propagation of shocks in this type of medium.
In regard to question (i), we present a range of experimental results
but only a partial answer, valid in a more restricted class of media.  In the
conclusion we discuss a potential research path to obtain such a condition for
more general media. 

We remark that linear wave propagation in a layered medium like that of
Figure \ref{fig:layered_intro} is well studied.  Due to the layering,
the leading-order effective speed of propagation depends on the angle; see e.g.
\cite{postma1955wave}.  Regarding the dispersion of linear waves in layered media,
see e.g. \cite{sun1968time,sun1968continuum,santosa1991,quezada2014two,farmer2018wave}.

The rest of this manuscript is organized as follows.
In \S\ref{sec:equations_and_material} we present the
system of equations and the type of media that we consider.
In \S\ref{sec:viscous_vs_disp_shocks} we discuss differences between
viscous and dispersive shocks. In \S\ref{sec:z-dispersion},
\ref{sec:c-dispersion} and \ref{sec:intro_shock_speed} we review
the different types of dispersive effects
that are induced in periodic media and, for a particular case, discuss
the speed of propagation of shocks proposed in \cite{Ketcheson_LeVeque_2011}. 
In \S\ref{sec:homogenized_equations} we use homogenization
theory to derive a leading order constant-coefficient system
that captures the main macroscopic features of the solution.
In \S\ref{sec:effective_RH} we propose an effective shock speed
by applying the Rankine-Hugoniot conditions to the homogenized system.
This estimate agrees with \cite{Ketcheson_LeVeque_2011} when $\theta=90^{\circ}$.
We test this estimate via numerical experiments, considering two different nonlinearities
and two types of periodic media. 
In \S\ref{sec:lax_entropy} we discuss possible conditions for
a shock to propagate without being regularized by the induced dispersion.
Finally, in \S\ref{sec:conclusions} we present some conclusions and
open questions, with possible approaches to resolve them.

% ******************************************************************** %
% ********** GOVERNING EQUATIONS AND MATERIAL CONFIGURATION ********** %
% ******************************************************************** %
\subsection{Governing equations and material configuration}\label{sec:equations_and_material}
We consider the scalar, nonlinear, variable-coefficient wave equation
\begin{align}\label{p-system}
  \epsilon_{tt}-\nabla\cdot\left(\frac{1}{\rho(\bfx)} \nabla\sigma(\epsilon,\bfx)\right) & = 0.
\end{align}
We use the notation of elasticity for consistency with related work
\cite{leveque2003, Ketcheson_LeVeque_2011, quezada2014two, ketcheson2015diffractons}.
Therefore $\epsilon$ denotes strain, $\rho$ density, and $\sigma$ stress.
Furthermore, let $K(\bfx)$ denote the bulk-modulus. We consider nonlinear stress-strain relations of the form
\begin{align}\label{nonlinear_const_relations}
  \sigma(\epsilon,\bfx) = \sigmahat(K(\bfx)\epsilon),
\end{align}
such that $\sigma_\epsilon$ can be expressed as $\sigma_\epsilon=G(\sigma)K$ for some function $G(\sigma)$.
This condition is needed during the homogenization process in \S\ref{sec:homogenized_equations}.
Note that this can be achieved if
$\sigmahat\in\mathcal{C}^1(K\epsilon)$ (i.e., if $\sigmahat$ is once differentiable with respect to $K\epsilon$)
and $\sigmahat$ is one-to-one.
%In this case $G(\sigma)=\sigmahat^\prime(K\epsilon)=\sigmahat^\prime\left(\sigmahat^{-1}(\sigma)\right)$.
%

Solutions of \eqref{p-system} with nonlinear stress-strain relations may 
involve shock singularities. In order to determine entropy-satisfying weak solutions, we write \eqref{p-system} 
as a first-order hyperbolic system in conservation form:
\begin{subequations}\label{p-system_conservation_form}
  \begin{equation}
    {\bf q}_{t}+{\bf f}({\bf q}, \bfx)_x+{\bf g}({\bf q}, \bfx)_y={\bf 0},
  \end{equation}
where
  \begin{align}
    {\bf q} & = \begin{bmatrix} \epsilon\\ \rho(\bfx) u \\ \rho(\bfx)  v \end{bmatrix},
    & {\bf f}({\bf q}, \bfx) & =\begin{bmatrix} -u \\ -\sigma(\epsilon,\bfx)\\ 0\end{bmatrix},
    & {\bf g}({\bf q}, \bfx) & =\begin{bmatrix} -v\\ 0\\ -\sigma(\epsilon,\bfx)\end{bmatrix}.
  \end{align}
\end{subequations}
Here $u$ and $v$ are the $x$- and $y$-components of velocity,
${\bf q}$ is the vector of conserved quantities, and ${\bf f}, {\bf g}$
are the components of the flux in the $x$- and $y$-directions respectively.
Note that (for simplicity) we concentrate on two dimensional waves but the results
we present can be readily extended to three dimensions. 
  We remark that system \eqref{p-system_conservation_form} is a two-dimensional
  generalization of the $p$-system from Lagrangian gas dynamics (see e.g.
  \cite[Section~2.13]{levequefvmbook}).

The medium is periodic and extends infinitely in both coordinate directions, as
shown in Figure \ref{fig:layered_intro}.
We consider both layered and smoothly-varying media.
The smoothly-varying medium is given by
\begin{subequations}\label{sin_medium}
\begin{align}
  K(\bfx) &= \frac{K_A+K_B}{2}+\frac{K_A-K_B}{2}\sin(2\pi\xi(\bfx)), \\
  \rho(\bfx) &= \frac{\rho_A+\rho_B}{2}+\frac{\rho_A-\rho_B}{2}\sin(2\pi\xi(\bfx)),
\end{align}
\end{subequations}
where $K_A$, $\rho_A$, $K_B$ and $\rho_B$ are strictly positive constants
and $\xi(\bfx)=x\sin\theta+y\cos\theta$.
The layered medium, which is composed of two types of materials: $A$ and $B$, is given by
\begin{align}
\label{layered_medium}
  (K(\bfx),\rho(\bfx)) &= \begin{cases} (K_A, \rho_A) & \text{ if } \sin(2\pi\xi(\bfx)) \ge 0 \\
                        (K_B, \rho_B) & \text{ if } \sin(2\pi\xi(\bfx)) < 0. \\\end{cases}
\end{align}

%\begin{figure}
%\begin{centering}
%  \includegraphics[scale=0.75]{figures/layered.png}
%\par\end{centering}
%\caption{\small Wave propagation in a layered medium in two dimensions. 
%  The material shown is repeated periodically and extends infinitely in both directions.
%  The direction of the heterogeneity is controlled by the angle $\theta$.
%  If $\theta=90^{\circ}$ we refer to the propagation as being transverse to the material heterogeneity 
%  and if $\theta=0^{\circ}$ we refer to the propagation as being parallel. 
%  \label{fig:layered}}
%\end{figure}

Important properties of the medium are characterized by the linearized impedance
$Z(\bfx):=\sqrt{K\rho}$ and sound speed $c(\bfx):=\sqrt{K/\rho}$.
Variations in $Z$ govern reflection, while $c$ dictates the speed of
propagation of small amplitude waves.
As we explain in \S \ref{sec:z-dispersion} and \ref{sec:c-dispersion},
different sources of dispersion are induced
due to changes in the impedance and sound speed.
In some of the numerical experiments that we perform we define
$K_A$, $\rho_A$, $K_B$ and $\rho_B$ by selecting the values of
$Z_A:=\sqrt{K_A\rho_A}$, $c_A:=\sqrt{K_A/\rho_A}$,
$Z_B:=\sqrt{K_B\rho_B}$ and $c_B:=\sqrt{K_B/\rho_B}$.

\subsubsection{Normalization of the material parameters}
Since the stress-strain relation $\sigma=\sigmahat(K\epsilon)$ is a function of the product $K\epsilon$,
we can normalize the parameters with respect to $(K_A,\rho_A)$.
As a result, we can consider without loss of generality $(K_A,\rho_A)=1$.
To see this, let 
\begin{subequations}\label{K_and_rho_tilde}
  \begin{align}
    \tilde{K}(\bfx)    &= K(\bfx)/K_A, \\
    \tilde{\rho}(\bfx) &= \rho(\bfx)/\rho_A.
    %  \begin{cases} 1 & \mbox{ } \forall \bfx\in\Omega_A\subset\Omega, \\ 
    %    \rho_B/\rho_A & \mbox{ } \forall \bfx\in\Omega_B\subset\Omega, 
    %  \end{cases}
  \end{align}
\end{subequations}
Multiply \eqref{p-system} by $K_A$ and use \eqref{K_and_rho_tilde} to obtain 
\begin{align}
  (K_A\epsilon)_{tt}-\nabla\cdot\left[\frac{1}{\tilde{\rho}(\bfx)} \nabla 
    \left( \frac{K_A}{\rho_A}\hat{\sigma}\left(\tilde{K}(\bfx)K_A\epsilon\right)\right)\right] & = 0.
\end{align}
Finally, let $\tilde{\epsilon}=K_A\epsilon$ and
$\tilde{\sigma}=\frac{K_A}{\rho_A}\hat{\sigma}\left(\tilde{K}(\bfx)K_A\epsilon\right)=\frac{K_A}{\rho_A}\sigma$
to obtain
\begin{align}
  \tilde{\epsilon}_{tt}-
  \nabla\cdot\left[\frac{1}{\tilde{\rho}(\bfx)} \nabla \tilde{\sigma}\left(\tilde{K}(\bfx)\tilde{\epsilon}\right)\right] = 0.
\end{align}
This equation has the same form as \eqref{p-system}, but with $\tilde{\rho}_A=\tilde{K}_A=1$
and the parameters in the $B$-material scaled accordingly.
Henceforth, we omit the tildes and assume $(K_A, \rho_A)=1$ whenever convenient.

% ****************************************************** %
% ********** VISCOUS VERSUS DISPERSIVE SHOCKS ********** %
% ****************************************************** %
\subsection{Viscous versus dispersive shocks}\label{sec:viscous_vs_disp_shocks}
In contrast to first-order hyperbolic conservation laws, higher-order PDE
models often include viscous or dispersive terms that prevent the formation of
discontinuities and lead instead to what are known as {\em viscous shocks} or
{\em dispersive shocks}.
Viscous and dispersive shocks have very different structures and are typically
treated with different mathematical tools.
Viscous shocks connect two states via a narrow, smooth transition region
and generate an increase in physical entropy over time.
The speed of a viscous shock is given by the Rankine-Hugoniot condition which
is a consequence of the conservation laws themselves \cite{lax1973hyperbolic}.
A dispersive shock, on the other hand, connects two states by an oscillatory
region that expands in time. Dispersive shocks do not generate entropy.
The speed of the leading edge of a dispersive shock is related to the speed of
a solitary wave, while the speed of the trailing edge is related to that of small-amplitude
perturbations.
The regularization observed in
\cite{santosa1991,leveque2003,Ketcheson_LeVeque_2011,quezada2014two,ketcheson2015diffractons}
and in the present
work can be related to dispersive regularization occurring in higher-order
dispersive wave models.  In particular, as explained in \cite{santosa1991}, the
periodic structure
%periodic variation
of the medium leads to an effective dispersion of waves.
We refer to \cite{karpman2016non} for a detailed discussion about viscous and dispersive shocks,
to \cite{whitham1965non} for the main mathematical tools to treat dispersive shocks, and to 
\cite{el2016dispersive} for a review of the theory of dispersive shocks. 
We remark that randomly varying media seem to lead to viscous, rather than
dispersive, effective equations; see for example \cite{solna2000ray,fouque2004shock,garnier2007effective}.

Numerical methods for first-order hyperbolic conservation laws are designed to
approximate the vanishing-viscosity limit of a viscously-regularized system.  Thus
they impose the Rankine-Hugoniot conditions and attempt to approximate viscous shocks
in the limit where the width of the shock region vanishes (i.e., the viscous shock
becomes a discontinuity).  
%We will refer to this limiting discontinuity as a viscous
%shock, even though it appears in the zero-viscosity limit, in order to distinguish it
%from a dispersive shock.  
In practice, any stable numerical method must impose some
viscosity in the shock region and so the shock is approximated by a region of finite
width.  Hence the solutions presented herein
are affected by an interplay of viscous (numerical) and dispersive (from the
medium) effects.  By refining the computational mesh, numerical viscosity can be reduced
and thus it is possible in principle to determine whether the solution to the first-order
variable-coefficient system (with no viscosity) is in fact NBAT.
We follow the methodology of \cite{Ketcheson_LeVeque_2011}, distinguishing
non-dispersive shocks by the presence of an increase in physical entropy
that persists when the mesh is refined.
%We estimate the speed of propagation of the viscous shock by applying the Rankine-Hugoniot
%conditions to a leading order homogenized constant-coefficient nonlinear system.

\begin{figure}
 \begin{centering}
   \includegraphics[scale=0.34]{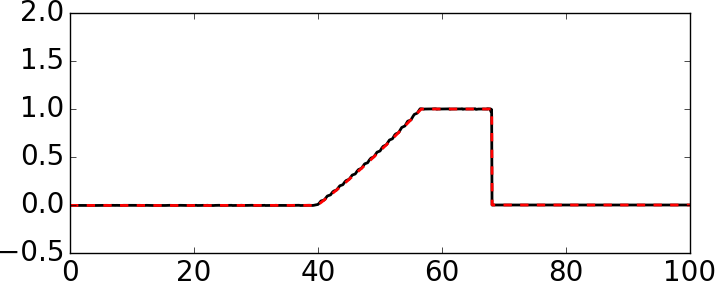}\quad
   \includegraphics[scale=0.34]{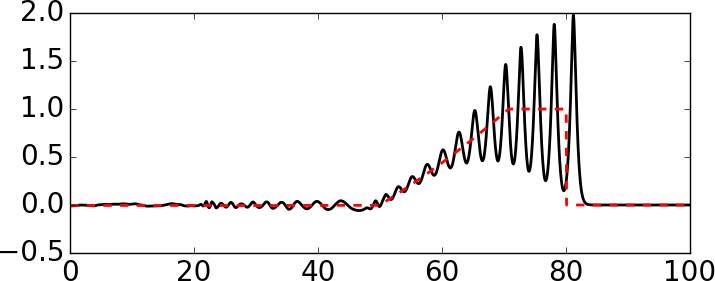}
   
   \includegraphics[scale=0.34]{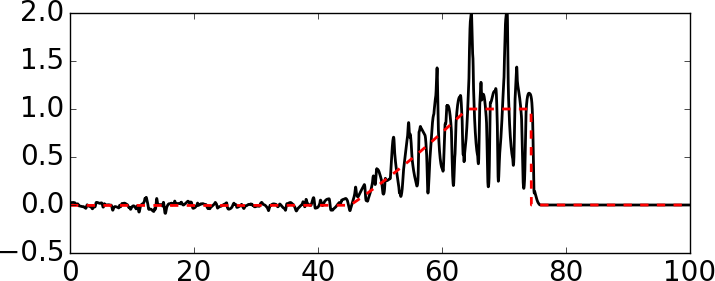}\quad
   \includegraphics[scale=0.34]{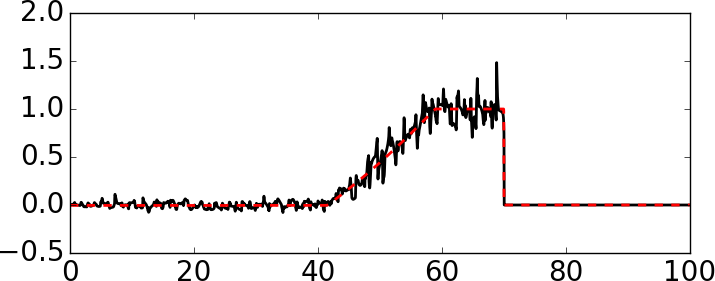}
  \par\end{centering}
\caption{\small Examples of wave propagation in a layered medium \eqref{layered_medium}.  Differences in material
  properties and angle of propagation can lead to behavior that is dominated by
  shock formation or by dispersion.
    The material coefficients are  $c_A=1$, $c_B=4$ and $Z_A=Z_B=1$
    and the angle of propagation is:
    (top-left) $\theta=90^\circ$,
    (top-right) $\theta=0^\circ$,
    (bottom-left) $\theta=45^\circ$ and
    (bottom-right) $\theta=67.5^\circ$.
    The dashed red line shows, for reference, the behavior
    in a homogeneous medium with appropriately averaged properties (see \S\ref{sec:homogenized_equations}).
    The initial condition is an effective purely right-going shock
    given by \eqref{init_cond_shock} with
    $x_s=20$, $u_l$ given by \eqref{right-going_shock}, $u_r=0$,
    $\sigma_l=1$ and $\sigma_r=0$.
    We show the solution at $t=25$ and use the nonlinear stress-strain relation \eqref{nonlinearity_exponential}.
  \label{fig:examples}}
\end{figure}

Figure \ref{fig:examples} demonstrates the range of possible behaviors.  Each plot
shows a solution corresponding to propagation of an initial square (plane) wave
in a periodic layered medium;
the only difference between the plots is $\theta$, the orientation
of the plane wave relative to the medium.  Depending on the medium properties and angle
of propagation, one may observe shock formation (top left), dispersive shock formation
followed by the creation of solitary waves (top right) or more complicated behavior
that lies between these extremes (bottom row).

% ********************************** %
% ********** Z-dispersion ********** %
% ********************************** %
\subsection{$Z$-dispersion: effective dispersion induced by reflections}\label{sec:z-dispersion}
\begin{figure}
\begin{centering}
  \includegraphics[scale=0.34]{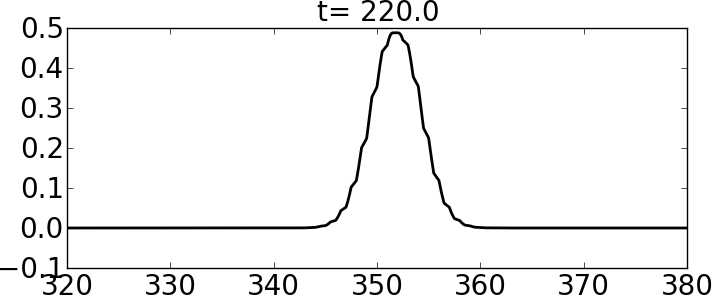}\quad
  \includegraphics[scale=0.34]{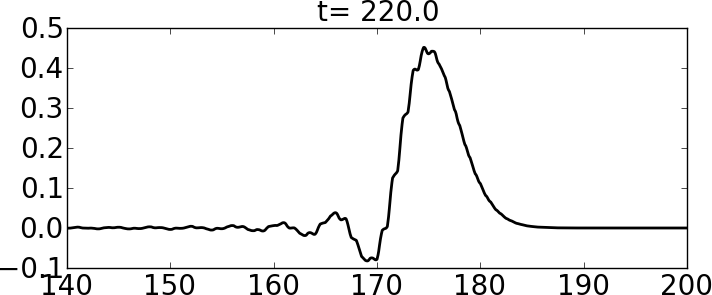}  
\par\end{centering}
\caption{\small Transverse ($\theta=90^\circ$) linear wave propagation in a layered medium \eqref{layered_medium}.
    The material coefficients are:
    (left) $c_A=1$, $c_B=4$ and $Z_A=Z_B=1$;
    (right) $c_A=c_B=1$, $Z_A=1$ and $Z_B=4$.
    The initial condition for both experiments is
    $\sigma(\bfx,t=0)=\exp(-x^2/10)$ and $u(\bfx,t=0)=v(\bfx,t=0)=0$.
  \label{fig:linear_Z-dispersive_media}}
\end{figure}

\begin{figure}
\begin{centering}
   	\includegraphics[scale=0.34]{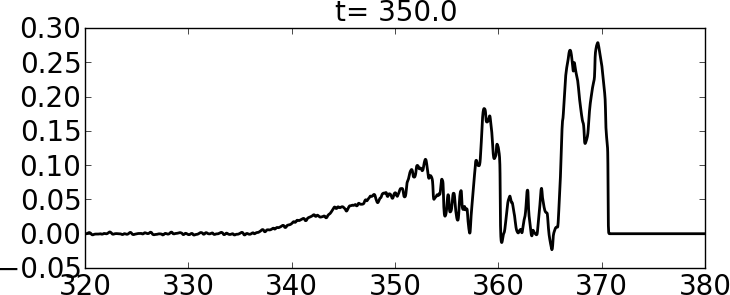}\quad
  	\includegraphics[scale=0.34]{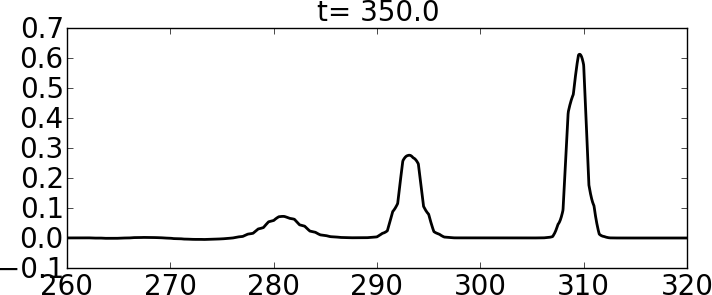}
\par\end{centering}
\caption{\small Transverse ($\theta=90^\circ$) nonlinear wave propagation in a layered medium \eqref{layered_medium}.
    The material coefficients are:
    (left) $c_A=c_B=1$, $Z_A=1$ and $Z_B=1.5$;
    (right) $c_A=c_B=1$, $Z_A=1$ and $Z_B=4$.
    The initial condition for both experiments is
    $\sigma(\bfx,t=0)=\exp(-x^2/10)$ and $u(\bfx,t=0)=v(\bfx,t=0)=0$.
  \label{fig:nonlinear_on_Z-dispersive_media}}
\end{figure}

Consider a plane wave traveling transversely in a layered medium; i.e. let
$\theta=90^{\circ}$ in Figure \ref{fig:layered_intro}.
In this case its propagation follows the theory developed in \cite{Ketcheson_LeVeque_2011},
which we briefly review here.
If the impedance is mismatched, reflections occur at each interface.
The net effect of these reflections is macroscopic dispersion, as shown in \cite{santosa1991}
via Bloch expansions and in \cite{Fish2001, conca1997homogenization, chen2001dispersive, Yong} via
homogenization theory. 
When the material parameters change such that the linearized impedance remains constant,
there are almost no reflections and, as a consequence, the dispersion introduced is negligible.
Moreover, if the material response is linear, the dispersion vanishes completely.
Therefore, having variable impedance is crucial to obtain effective dispersion
in waves traveling transverse to the material heterogeneity. 
We show this in the left and right panels of Figure \ref{fig:linear_Z-dispersive_media} by considering
a linear wave traveling in a layered medium with constant and mismatched impedance, respectively.
With nonlinear waves, if the variations in impedance are small, shocks may develop in finite time
(like in a homogeneous medium). This is shown in the left panel of Figure
\ref{fig:nonlinear_on_Z-dispersive_media}.
If the impedance mismatch is large, the effective dispersion may be strong
enough to delay or avoid shock formation \cite{Ketcheson_LeVeque_2011}; i.e.,
the induced dispersion can act as a regularization mechanism.
Moreover, the nonlinear and the dispersive effects can balance each other leading to solitary
wave formation \cite{leveque2003}, as seen in the right panel of Figure \ref{fig:nonlinear_on_Z-dispersive_media}.
This is similar to the formation of solitons in the KdV equation, see for instance \cite{zabusky1965interaction}.
However, it is important to remember that in the present setting the PDE being solved
contains no dispersive terms.
%
%In \cite{2014_cylindrical_solitary_waves} the authors used finite volume methods to obtain cylindrical solitary waves
%on two dimensional periodic media arranged in a checkerboard pattern. The material parameters
%are selected such that the impedance variations but the sound speed remains constant. 
%Not surprisingly, if the impedance contrast is not large enough, shocks start to appear;
%otherwise, the initial condition evolves into solitary waves. We illustrate these situations in Figure
%\ref{fig:nonlinear_Z-dispersive_media_on_checker}.

%% \begin{figure}
%% \begin{centering}
%%  \subfloat[Initial condition\label{fig:2D_stegotons_init_condition}]{
%%  	\includegraphics[scale=0.2]{figures/2D_stegotons_initial_condition_frame_0.png}}
%%  \subfloat[Checkerboard domain\label{fig:checker}]{
%%  	\includegraphics[scale=0.2]{figures/checker.png}}
%%  \subfloat[Two-dimensional nonlinear wave in a checker board domain. \label{fig:2D_stegoton_and_shock}]{
%%   	\includegraphics[scale=0.2]{figures/2D_shock_t40_frame_40.png}
%%     	\includegraphics[scale=0.2]{figures/2D_stegotons_t40_frame_40.png}}    
%% \par\end{centering}
%% \caption{\small Nonlinear wave propagation in a two-dimensional checkerboard medium with changing impedance (and constant sound speed). 
%%   (a) Zoomed initial solution.
%%   (b) Zoomed checkerboard medium. 
%%   (c) Solution at $t=40$ for a medium with impedance contrast of $1.5$ and $4$ in the left and right panels respectively.
%%   See \S \ref{sec:equations_and_material} for details about the equations and material parameters. 
%% \label{fig:nonlinear_Z-dispersive_media_on_checker}}
%% \end{figure}

Up to this point we have considered dispersion induced by effective reflections at the material interfaces when the
impedance varies. We refer to this effect as {\it `$Z$-dispersion'} since it depends on variation in the impedance,
which is commonly denoted by $Z$ in elasticity theory.

% ********************************** %
% ********** C-DISPERSION ********** %
% ********************************** %
\subsection{$C$-dispersion: effective dispersion induced by diffraction}\label{sec:c-dispersion}
\begin{figure}
  \begin{centering}
    \includegraphics[scale=0.33]{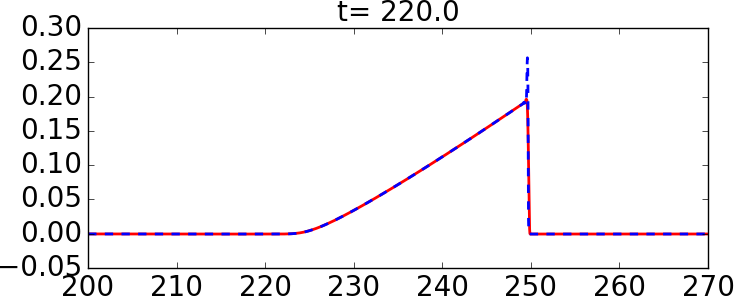}\qquad 
    \includegraphics[scale=0.33]{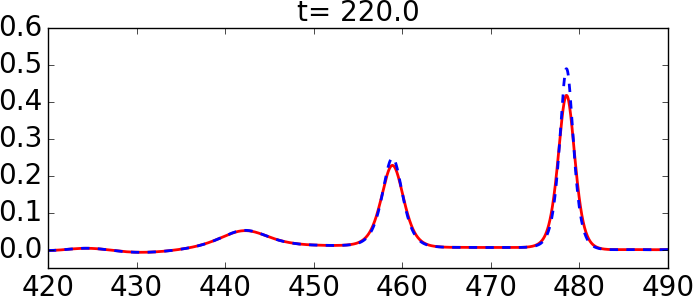}
    \par\end{centering}
    \caption{\small
      Parallel ($\theta=0^\circ$) nonlinear wave propagation in a layered medium \eqref{layered_medium}.
      Sufficiently large variation in sound speed seems to prevent shock formation.
        The material coefficients are:
        (left) $c_A=1$, $c_B=1.1$ and $Z_A=Z_B=1$;
        (right) $c_A=1$, $c_B=4$ and $Z_A=Z_B=1$.
          We plot in solid red a slice through the middle of the $B$-layer
          and in dashed blue a slice through the middle of the $A$-layer.
        The initial condition for both experiments is
        $\sigma(\bfx,t=0)=\exp(-x^2/10)$ and $u(\bfx,t=0)=v(\bfx,t=0)=0$.
      \label{fig:nonlinear_c-dispersive_media}}
\end{figure}

In \cite{quezada2014two} the authors considered linear waves traveling in the same type of medium
considered in the present work.  Let us now review what was observed therein for parallel propagation;
i.e. $\theta=0^{\circ}$ in Figure \ref{fig:layered_intro}.   In this case, if the sound speed differs
in the $A$ and $B$ layers, effective dispersion is again introduced (even if the impedance is constant).
We refer to this effect as {\it `$C$-dispersion'} since it depends on variation in the sound speed,
which is commonly denoted by $c$ in elasticity theory. 
In the nonlinear case, $C$-dispersion can also act as a regularization
mechanism to delay or avoid the formation of shocks. 
If the dispersion is large, the nonlinear and the dispersive effects can balance each other,
leading to the formation of solitary waves \cite{ketcheson2015diffractons}.
In Figure \ref{fig:nonlinear_c-dispersive_media} we demonstrate this effect.

%% \begin{figure}
%%   \begin{centering}
%%     \includegraphics[scale=0.4]{figures/het_matchC_linear_stress_slices}
%%     \includegraphics[scale=0.4]{figures/het_mismatchC_linear_stress_slices}
%%     \par\end{centering}
%%     \caption{\small Linear waves traveling parallel to the material heterogeneity in a layered medium.
%%       In the left and right panels we consider constant and mismatched sound speed respectively. 
%%       The sound speed contrast in the right panel is of 4. 
%%       We plot slices at the middle of each material section in the layered medium. 
%%       See \S \ref{sec:equations_and_material} for details about the equations and material parameters. 
%% \label{fig:linear_c-dispersive_media}}
%% \end{figure}

For general propagation along an arbitrary angle $\theta$, effective $Z$-dispersion 
and effective $C$-dispersion may both be introduced,
due to variations in the impedance
and the sound speed, respectively.  The strength of each effect depends on the material
parameters and the direction of propagation.
\subsection{Shock speed}\label{sec:intro_shock_speed}
\begin{figure}
 \begin{centering}
   \subfloat[Layered medium with $\theta=90^\circ$\label{fig:right_going_shock_theta90_intro}]{
     \includegraphics[scale=0.34]{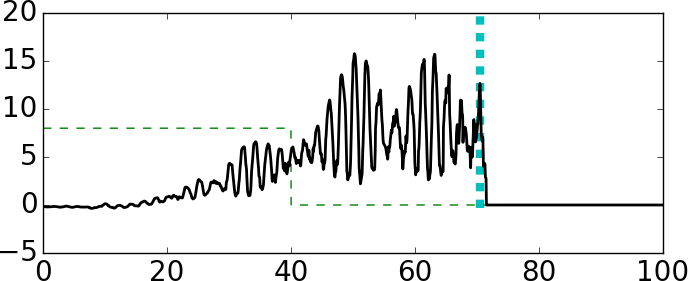}}
   \quad
   \subfloat[Layered medium with $\theta=45^\circ$\label{fig:right_going_shock_theta45_intro}]{
     \includegraphics[scale=0.34]{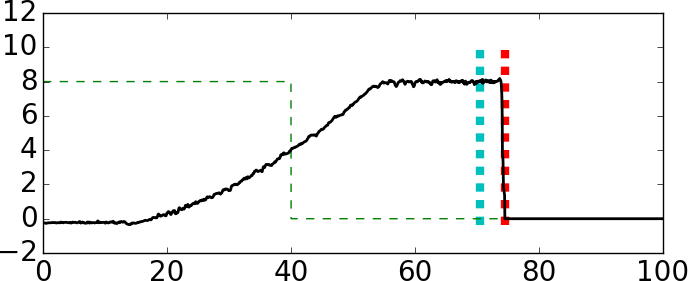}}
  \par\end{centering}
  \caption{\small Right-going shock in a layered medium \eqref{layered_medium} with 
      (a) $\theta=90^\circ$ and (b) $\theta=45^\circ$.
      The material coefficients are $Z_B=1$, $Z_A=4$ and $c_A=c_B=1$.    
      The initial condition is an effective purely right-going shock
      given by \eqref{init_cond_shock} with
      $x_s=40$, $u_l$ given by \eqref{right-going_shock}, $u_r=0$,
      $\sigma_l=8$ and $\sigma_r=0$.
      We show the solution at $t=20$ and use the nonlinear stress-strain relation \eqref{nonlinearity_exponential}
    In all plots, the black solid line is the solution along $y=0.5$,
    the cyan dashed line is the shock position predicted by \eqref{speed_of_shocks_z_dispersion} and
    the green dashed line is the initial condition. 
    In (b), the red dashed line is the shock position predicted by our estimate
    \eqref{speed_of_shock}, which we present in \S\ref{sec:effective_RH}.
    \label{fig:right-going_shocks_intro}}
\end{figure}

In this work we focus on a setting in which the dispersive effects are
small (or the initial data is large), so that shocks may still appear. We are
interested in the speed of propagation of the resulting shocks. 
This was studied in \cite{Ketcheson_LeVeque_2011} for the case of transverse
propagation (i.e., $\theta=90^{\circ}$ in Figure \ref{fig:layered_intro}),
in which case the problem is one-dimensional and there are only
$Z$-dispersive effects.
Following the notation therein, let $\sigma$, $\epsilon$ and $\rho$
denote stress, strain, and density, respectively. The hypothesized shock
speed from \cite{Ketcheson_LeVeque_2011} is then
\begin{align}\label{speed_of_shocks_z_dispersion}
  s=\sqrt{\frac{1}{\rho_m}\left(\frac{[\sigma]}{[\epsilon]}\right)_h},
\end{align}
where $(\cdot)_h$ is the harmonic average operator (over one spatial period of the medium),
$\rho_m$ is the mean density and $[z]=z_r-z_l$ denotes the jump in quantity $z$ across the shock.
Formula \eqref{speed_of_shocks_z_dispersion} is just the usual Rankine-Hugoniot shock speed
but with each quantity replaced by an appropriate spatial average.  The choice to use
an ordinary average for $\rho$ and a harmonic average for $\sigma/\epsilon$ is based on
the fact that small-amplitude, long-wavelength pulses travel at an effective speed given by 
\begin{align}
  c_{\text{eff}}=\sqrt{\frac{K_h}{\rho_m}},
\end{align}
where $K$ is the bulk-modulus; see \cite{santosa1991}.
Thus \eqref{speed_of_shocks_z_dispersion} is a natural extension of the Rankine-Hugoniot
condition to determine the speed of shocks using effective material parameters.
In Figure \ref{fig:right_going_shock_theta90_intro} we demonstrate the correctness of this estimate by considering
a right going shock (as those we study in \S \ref{sec:right-going_shock}) propagating in a 
$Z$-dispersive medium with impedance in material A and B to be given by
$Z_A=1$ and $Z_B=4$ respectively. The cyan dashed line represents the position of the shock
as predicted by \eqref{speed_of_shocks_z_dispersion}.
%
%We remark from \eqref{speed_of_shocks_z_dispersion} that the speed of the shock cannot
%be given by a simple average of the shock speeds in the different sections of the layered medium. 

One might expect that, under the more general scenario in Figure \ref{fig:layered_intro},
shocks would propagate with the same speed \eqref{speed_of_shocks_z_dispersion}.
However, one quickly finds that this is not the case.
For example, in Figure \ref{fig:right_going_shock_theta45_intro} we consider
propagation along $\theta=45^{\circ}$ in a medium with the same material properties as before.
Again, we indicate the shock position predicted by \eqref{speed_of_shocks_z_dispersion} with
a cyan dashed line. It is clear that the results in \cite{Ketcheson_LeVeque_2011}
do not apply for the more general framework that we study in this work.
In Figure \ref{fig:right_going_shock_theta45_intro} we show
with a dashed red line the position predicted by the theory developed in the present work.
We propose an effective speed of propagation that is based on a leading order homogenized system.

% ******************************************* %
% ********** HOMOGENIZED EQUATIONS ********** %
% ******************************************* %
\section{Homogenized equations}\label{sec:homogenized_equations}
Our first goal is to generalize the proposed shock speed that was given in
\cite{Ketcheson_LeVeque_2011} for the case $\theta=90^\circ$.
In order to do so, and also to place that hypothesis on a firmer
mathematical footing, we perform an asymptotic analysis.
Specifically, in this section we use homogenization theory to capture the
macroscopic effects when waves with large
wavelength travel in a
periodic medium with small period $\Omega$ (relative to the wavelength).
The leading-order result of this process is a constant-coefficient system
with the material properties given by some average that depends on the relative
angle between the propagation of the wave and the variation in the medium.
We follow \cite{ketcheson2015diffractons} and references therein. 
Recall system \eqref{p-system_conservation_form}:
%Consider \eqref{p-system} as a first order system of conservation laws: 
\begin{subequations}\label{p-system_cons_form_hom_section}
\begin{align}
	\epsilon_t-(u_x+v_y) & = 0, \\
	\rho u_t - \sigma_x & = 0, \\
	\rho v_t - \sigma_y & = 0.
\end{align}
\end{subequations}
From \S\ref{sec:equations_and_material} we have $\sigma_\epsilon=K(\bfx)G(\sigma)$.
Using the chain rule, $\sigma_t=\sigma_\epsilon\epsilon_t$. Therefore, we can write \eqref{p-system_cons_form_hom_section} as
\begin{subequations}\label{system_to_homogenized}
\begin{align}
	K^{-1}\sigma_t-G(\sigma)(u_x+v_y) & = 0, \\
	\rho u_t - \sigma_x & = 0, \\
	\rho v_t - \sigma_y & = 0.
\end{align}
\end{subequations}
We start by introducing a parameter $\delta=\Omega/\lambda$, where
$\Omega$ is the material period and $\lambda$ is the wave length.
We restrict our attention to waves with $\lambda\gg\Omega\implies\delta\ll 1$.
We recognize and introduce a fast spatial scale
$\hat{\xi}=\delta^{-1}\xi$ where $\xi=x\sin\theta+y\cos\theta$
(is obtained by a simple rotation of axes and) defines the direction of the heterogeneity. 
We assume $\sigma=\sigma(x,y,\hat{\xi},t)$, $u=u(x,y,\hat{\xi},t)$ and $v=v(x,y,\hat{\xi},t)$
and that the material properties depend only on the fast scale; i.e., 
$\rho=\rho(\hat{\xi})$ and $K=K(\hat{\xi})$.

By the chain rule $\partial_x\mapsto\partial_x+\delta^{-1}\partial_{\xihat}\sin\theta$
and $\partial_y\mapsto\partial_y+\delta^{-1}\partial_{\xihat}\cos\theta$.
Therefore \eqref{system_to_homogenized} becomes
\begin{subequations}\label{scaled_system}
  \begin{align}
    K^{-1}\sigma_t-G(\sigma)[u_x+v_y+\delta^{-1}(u_{\xihat}\sin\theta+v_{\xihat}\cos\theta)] & = 0, \\
    \rho u_t - (\sigma_x+\delta^{-1}\sigma_{\xihat}\sin\theta) & = 0, \\
    \rho v_t - (\sigma_y+\delta^{-1}\sigma_{\xihat}\cos\theta) & = 0.
\end{align}
\end{subequations}
The next step is to expand $\sigma(x,y,\hat{\xi},t)$, $u(x,y,\hat{\xi},t)$ and $v(x,y,\hat{\xi},t)$ using
the small parameter $\delta$. For example,
$\sigma(x,y,\hat{\xi},t)=\sum_{i=0}^\infty \delta^i \sigma_i (x,y,\hat{\xi},t)$ 
and similarly for $u$ and $v$. We plug these expansions into \eqref{scaled_system} to get
\begin{subequations} \label{expanded_system}
\begin{align}
 K^{-1}\sum_{i=0}^\infty\delta^i\sigma_{i,t}
 -G\left(\sigma\right)\left(\sum_{i=0}^\infty\delta^i u_{i,x}+\sum_{i=0}^\infty\delta^i v_{i,y}
 +\sum_{i=0}^\infty \delta^{i-1}(u_{i,\xihat}\sin\theta+v_{i,\xihat}\cos\theta) \right) & = \
0, \\
 \rho\sum_{i=0}^\infty\delta^i u_{i,t}-\left(\sum_{i=0}^\infty\delta^i \sigma_{i,x}+\sum_{i=0}^\infty \delta^{i-1}\sigma_{i,\xihat}\sin\theta\right) & = 0, \\
 \rho\sum_{i=0}^\infty\delta^i v_{i,t}-\left(\sum_{i=0}^\infty\delta^i \sigma_{i,y}+\sum_{i=0}^\infty \delta^{i-1}\sigma_{i,\xihat}\cos\theta\right) & = 0,
\end{align}
\end{subequations}
where $(\cdot)_{i,z}$ denotes differentiation of $(\cdot)_i$ with respect to $z$.
The function $G(\sigma)$ is expanded around $\sigma_0$ using Taylor series as
$G(\sigma)=G(\sigma_0)+\delta G'(\sigma_0) \sigma_1 + \dots$.

From \eqref{expanded_system}, we collect the terms of order $\order(\delta^{-1})$:
\begin{subequations} \label{order delta^-1 system}
\begin{align}
	(u_0\sin\theta+v_0\cos\theta)_{\hat{\xi}} & = 0, \label{order_delta^-1_system_eqn1} \\
	\sigma_{0,\hat{\xi}}\sin\theta & = 0, \\
	\sigma_{0,\hat{\xi}}\cos\theta & = 0,
\end{align}
\end{subequations}
which implies that $\sigma_0=:\bar{\sigma}_0(x,y,t)$
and that $u_0\sin\theta+v_0\cos\theta$ is independent of $\xihat$.
Here we use the notation that variables independent of the fast scale $\xihat$ are denoted by a bar. 

From \eqref{expanded_system}, we collect the terms of order $\order(\delta)$: 
\begin{subequations} \label{order 1 system}
\begin{align}
  K^{-1}\bar{\sigma}_{0,t}-G(\bar{\sigma}_0)\left(u_{0,x}+v_{0,y}+u_{1,\xihat}\sin\theta+v_{1,\xihat}\cos\theta \right) &=0,\label{order_1_system_eqn1}\\
  \rho u_{0,t}-\left(\bar{\sigma}_{0,x}+\sigma_{1,\xihat}\sin\theta\right) &= 0, \label{order_1_system_eqn2}\\
  \rho v_{0,t}-\left(\bar{\sigma}_{0,y}+\sigma_{1,\xihat}\cos\theta\right) &= 0. \label{order_1_system_eqn3}
\end{align}
\end{subequations}
Multiply equation \eqref{order_1_system_eqn2} by $\rho^{-1}$, apply the averaging operator 
\begin{align}
  \avg[\cdot]:=\frac{1}{\Omega}\int_0^{\Omega}(\cdot)d\xi
  %=\frac{1}{\lambda\sin\theta}\int_0^{\lambda\sin\theta}(\cdot)d\xihat
  =\frac{1}{\lambda}\int_0^{\lambda}(\cdot)d\xihat
\end{align}
and multiply the resulting equation by $\cos\theta$.
Similarly, multiply equation \eqref{order_1_system_eqn3} by $\rho^{-1}$, 
apply the averaging operator and multiply the resulting equation by $\sin\theta$. 
Combine the two equations to get
\begin{align}\label{homog_aux1}
  \left(\bar{u}_{0}\cos\theta - \bar{v}_{0}\sin\theta \right)_t-\rho_h^{-1}\left(\bar{\sigma}_{0,x}\cos\theta-\bar{\sigma}_{0,y}\sin\theta \right) = 0,
\end{align}
where $\rho_h:=\avg[\rho^{-1}]^{-1}$ denotes the harmonic average. 
Here we defined $\bar{u}_0:=\avg[u_0]$ and $\bar{v}_0:=\avg[v_0]$.
Multiply equation \eqref{order_1_system_eqn2} by $\sin\theta$ and 
equation \eqref{order_1_system_eqn3} by $\cos\theta$ and sum the two equations to obtain
\begin{align}\label{hom_aux1}
  \rho\left(u_0\sin\theta+v_0\cos\theta\right)_t-(\bar{\sigma}_{0,x}\sin\theta+\bar{\sigma}_{0,y}\cos\theta)-\sigma_{1,\hat{\xi}}=0.
\end{align}
Apply the averaging operator to \eqref{hom_aux1} noting that by \eqref{order_delta^-1_system_eqn1} 
$u_0\sin\theta+v_0\cos\theta$ is independent of $\hat{\xi}$.
Doing this yields
\begin{align} \label{homog_aux2}
  \left(\bar{u}_{0}\sin\theta + \bar{v}_{0}\cos\theta \right)_t-\rho_m^{-1}\left(\bar{\sigma}_{0,x}\sin\theta+\bar{\sigma}_{0,y}\cos\theta \right) = 0,
\end{align}
where $\rho_m:=\avg[\rho]$ denotes the arithmetic average. Here we assumed $\sigma_{1,\xihat}$ is periodic
with respect to the fast variable $\xihat$, which is a standard assumption in homogenization theory
(see for instance \cite{Fish2001, chen2001dispersive} and references therein). 
Multiply \eqref{homog_aux1} by $\cos\theta$ and \eqref{homog_aux2} by $\sin\theta$ 
and sum the two equations to obtain
\begin{align}
  \bar{u}_{0,t}-\left(\rho_h^{-1}\cos^2\theta+\rho_m^{-1}\sin^2\theta \right)\bar{\sigma}_{0,x}
  -\left(\rho_m^{-1}-\rho_h^{-1}\right)\sin\theta\cos\theta\bar{\sigma}_{0,y}
  = 0.
\end{align}
Apply the averaging operator to \eqref{order_1_system_eqn1}
(again assuming $\xihat$-periodicity of $u_{1,\xihat}$ and $v_{1,\xihat}$) to obtain
\begin{align}
  \Kbar^{-1}\bar{\sigma}_{0,t}-G(\bar{\sigma}_0)(\bar{u}_{0,x} + \bar{v}_{0,y}) &= 0.
\end{align}

Finally, for a plane wave propagating in the $x$-direction, we have $\bar{\sigma}_{0,y}, \bar{v}_{0,y}=0$.
Therefore, 
\begin{subequations}\label{homogenized_leading_order}
\begin{align}
  \Kbar^{-1}\bar{\sigma}_{0,t}-G(\bar{\sigma}_0)\bar{u}_{0,x} &= 0, \\
  \rhostar\bar{u}_{0,t}-\bar{\sigma}_{0,x} &= 0,
\end{align}
\end{subequations}
where $\rhostar=\left(\rho_h^{-1}\cos^2\theta+\rho_m^{-1}\sin^2\theta \right)^{-1}$ 
is the effective density in the oblique medium and 
$\Kbar=\avg[K^{-1}]^{-1}$ is the effective bulk-modulus.
Equations \eqref{homogenized_leading_order} constitute the leading order constant-coefficient homogenized system. 
We can now identify the effective material properties. The effective bulk-modulus is given 
by the harmonic average and the effective density is given by $\rhostar$.
This agrees with the theory in \cite{quezada2014two} where the two cases $\theta=0^\circ$ and $\theta=90^\circ$ are considered.

% ************************************************ %
% ********** EFFECTIVE RANKINE-HUGONIOT ********** %
% ************************************************ %
\section{Effective Rankine-Hugoniot conditions}\label{sec:effective_RH}
The homogenized equations, which are given by \eqref{homogenized_leading_order},
can be rewritten as a system of conservation laws: 
\begin{subequations}\label{homogenized_leading_order_in_conservation_form}
\begin{align}
  \epsilon_t-u_x & = 0, \\
  \rhostar u_t - \sigmabar_x & = 0.
\end{align}
\end{subequations}
The system \eqref{homogenized_leading_order_in_conservation_form} is closed by
choosing a nonlinear stress-strain relation $\sigmabar(\epsilon)=\sigma(\epsilon; \Kbar)$.
This leading order system is a dispersionless approximation of the original problem; i.e., it captures
the effective behavior of the solution neglecting the dispersive effects.
Therefore, it yields a better approximation when the dispersion is small.
In this work we are interested in the case when (viscous) shocks form which is expected
to happen when the effective dispersion is small.
In this regime, it is reasonable to expect the shock to propagate with speed close to
that given by applying the Rankine-Hugoniot conditions to \eqref{homogenized_leading_order_in_conservation_form}.
By doing this we obtain 
\begin{subequations}\label{effective_RH_conditions}
\begin{align}
  s[\epsilon] &=-[u], \label{effective_RH_conditions_eqn1}\\
  \rhostar s[u]  &=-[\sigmabar],
\end{align}
\end{subequations}
where $s$ is the speed of the shock and $[(\cdot)]$ is the jump of $(\cdot)$; 
i.e., $[(\cdot)]:=(\cdot)_l-(\cdot)_r$. Therefore, we propose an estimate
for the speed of propagation of viscous shocks to be given by 
\begin{align}\label{speed_of_shock}
  s_{\text{eff}}=\pm\sqrt{\frac{[\sigmabar]}{\rhostar[\epsilon]}},
\end{align}
where the minus and plus sign correspond to the left- and right-going shocks, respectively.
  Note that $\rhostar$ increases as $\theta$ gets closer to $90^\circ$.
  Therefore, the speed of the shock is reduced as $\theta$ gets closer to $90^\circ$.
  A similar result for long-wavelength linear waves is well known; see
  e.g. \cite{postma1955wave,backus1962long}.

\subsection{Effective purely right-going shocks}\label{sec:right-going_shock}
The Rankine-Hugoniot conditions not only provide the speed of propagation of shocks
but also the proper way to connect left and right states of a shock. 
Without loss of generality we consider just right-going shocks,
which correspond to the positive speed $s_{\text{eff}}$ in \eqref{speed_of_shock}.
From \eqref{effective_RH_conditions}, we get 
\begin{align}\label{right-going_shock}
  [u]=-\sqrt{\frac{[\sigmabar][\epsilon]}{\rhostar}}.
\end{align}
This expression provides the connection between the left and right state for a right-going shock. 
We perform experiments on different layered media with coefficients given by \eqref{layered_medium}.
In these experiments the nonlinear constitutive relation $\sigmabar(\epsilon)$
is given by \eqref{nonlinearity_exponential} with $K=\Kbar$, and the initial
condition is
\begin{align}\label{init_cond_shock}
  \sigma(\bfx,t=0)=
  \begin{cases}
    \sigma_l, &\text{ if } x\leq x_s \\
    \sigma_r, &\text{ if } x>x_s
  \end{cases},
  &&
  u(\bfx,t=0)=
  \begin{cases}
    u_l, &\text{ if } x\leq x_s \\
    u_r, &\text{ if } x>x_s
  \end{cases},
  &&
  v(\bfx,t=0)=0,
\end{align}
where $x_s$ defines the position of the shock,
$u_l$ is given by \eqref{right-going_shock}, $\sigma_l=1$, $\sigma_r=0$ and $u_r=0$. 
For all the experiments, we solve the original variable coefficient system \eqref{p-system} and the leading order
homogenized system \eqref{homogenized_leading_order_in_conservation_form}.
We consider two stress relations:
\begin{subequations}
\begin{align}
  \sigma(\epsilon,\bfx) &= \exp(K(\bfx)\epsilon) - 1, \label{nonlinearity_exponential}\\
  \sigma(\epsilon,\bfx) & = \alpha K(\bfx)\epsilon+\beta(K(\bfx)\epsilon)^2+\gamma(K(\bfx)\epsilon)^3, \label{nonlinearity_polynomial}
\end{align}
\end{subequations}
with $\alpha, \beta, \gamma\in\mathbb{R}$.
Both are smooth, one-to-one functions; the first has been used before in \cite{leveque2003}
while the second is important in nonlinear optics \cite{boyd2003nonlinear}.
We take $\alpha=0.1, \beta=0$ and $\gamma=5$.

\begin{figure}
 \begin{centering}
   \includegraphics[scale=0.175]{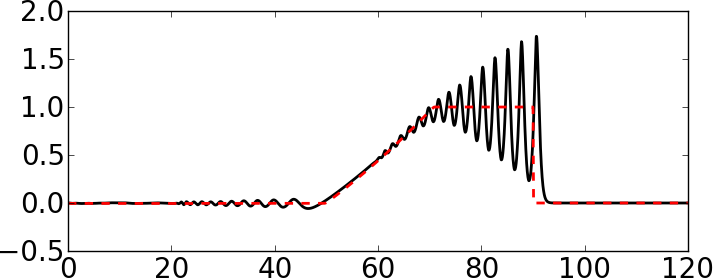}
   \includegraphics[scale=0.175]{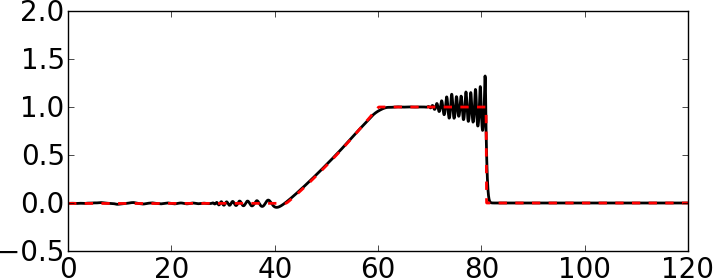}
   \begin{overpic}[scale=0.175]{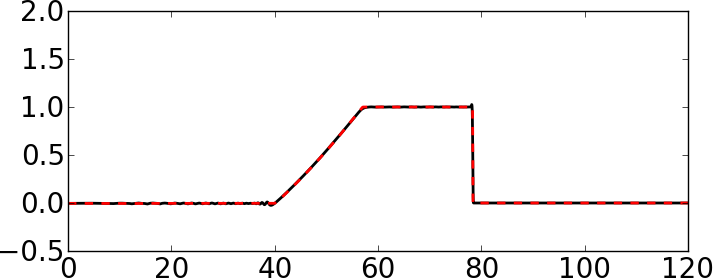}
     \put(52.5,28){\includegraphics[scale=0.15]{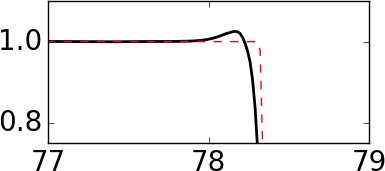}}
   \end{overpic}
   \begin{overpic}[scale=0.175]{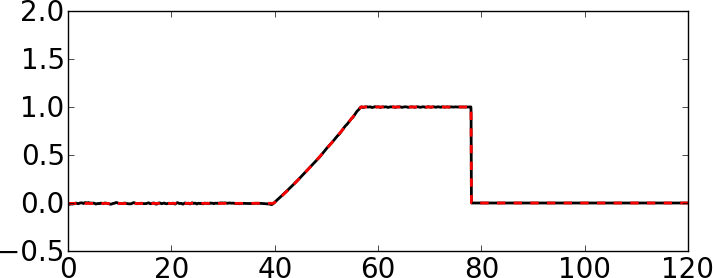}
     \put(52.5,28){\includegraphics[scale=0.15]{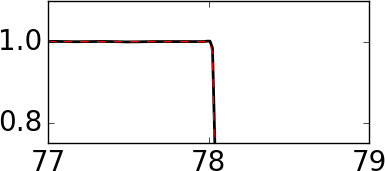}}
   \end{overpic}
  \par\end{centering}
  \caption{\small Effective right-going shocks in a layered medium with $\theta=0^\circ$.
      The material coefficients are given by \eqref{layered_medium} with
      sound speed contrast of (from left to right) $c_B/c_A=4, ~2, ~1.25$ and $1$.
      When $c_B/c_A=1$ we use $K_A=\rho_A=1$ and $K_B=\rho_B=4$. For all other cases
      we use \eqref{material_properties_rgoing_shocks_theta0}. 
      The red dashed line is the solution of the homogenized system \eqref{homogenized_leading_order}.    
      The initial condition is an effective purely right-going shock
      given by \eqref{init_cond_shock} with
      $x_s=30$, $u_l$ given by \eqref{right-going_shock}, $u_r=0$,
      $\sigma_l=1$ and $\sigma_r=0$.
      We show the solution at $t=40$ and use the nonlinear stress-strain relation \eqref{nonlinearity_exponential}.
      For the last two plots, we show an enlarged view of the region near the shock.
\label{fig:right-going_on_c-dispersive_media}}
\end{figure}

We first consider parallel propagation ($\theta=0^\circ$ in Figure \ref{fig:layered_intro}).
In this case $C$-dispersion is introduced by diffraction due to variation in the sound speed. 
Results are shown in Figure \ref{fig:right-going_on_c-dispersive_media} for
media with sound speed contrast of $c_B/c_A=4,~2,~1.25, \mbox{ and } 1$. 
When $c_B/c_A=1$ we use $K_A=\rho_A=1$ and $K_B=\rho_B=4$; i.e., 
the medium is still non-homogeneous. For all other cases we use
\begin{subequations}\label{material_properties_rgoing_shocks_theta0}
\begin{align}
  K_A &= \frac{1+c_B/c_A}{2c_B/c_A},  & \qquad \rho_A &= \frac{1}{K_A}, \\
  K_B &= \frac{1+c_B/c_A}{2},  & \qquad \rho_B &= \frac{1}{K_B}.
\end{align}
\end{subequations}
In the last three cases the approximation is very good.  In the first case,
a train of solitons forms and no singular shock is discernible.
Clearly, the approximation is better when the sound speed contrast is small.
For the last two plots, we show a zoomed image of the region near the shock.
Note that the effect of dispersion is still evident when $c_B/c_A=1.25$,
which is not the case when $c_B/c_A=1$.
  
\begin{figure}
 \begin{centering}
   \subfloat[$\theta=22.5^\circ$]{
     \includegraphics[scale=0.23]{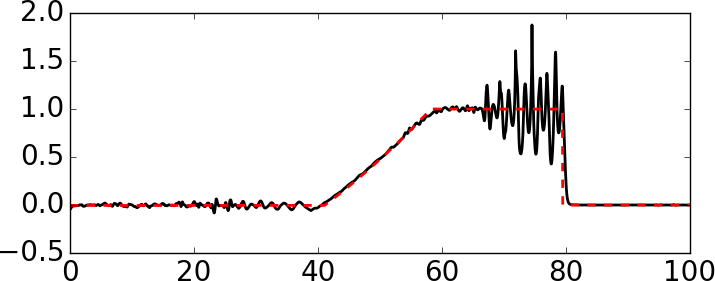}
     \includegraphics[scale=0.23]{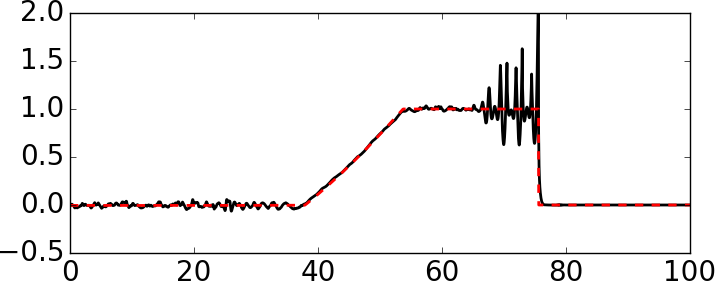}                
     \includegraphics[scale=0.23]{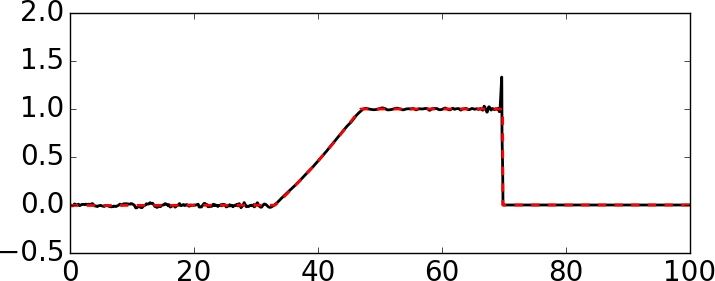}}                
   
   \subfloat[$\theta=45^\circ$]{
     \includegraphics[scale=0.23]{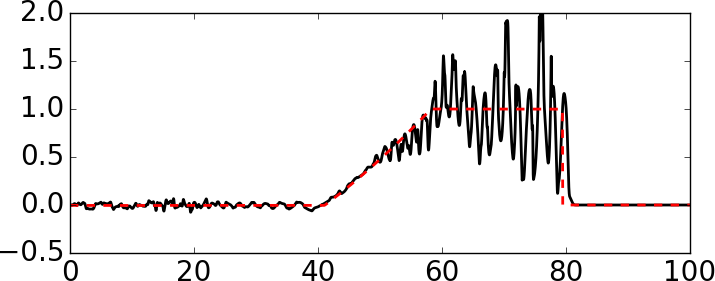}
     \includegraphics[scale=0.23]{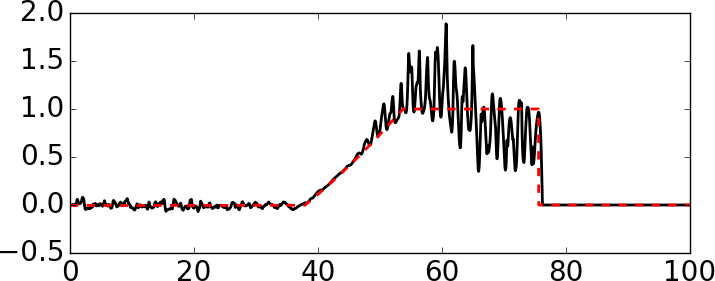}
     \includegraphics[scale=0.23]{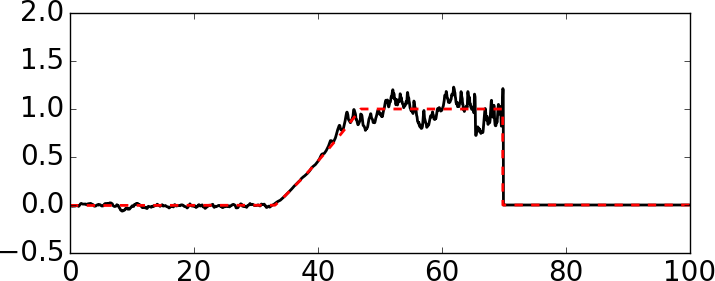}}
   
   \subfloat[$\theta=67.5^\circ$]{
     \includegraphics[scale=0.23]{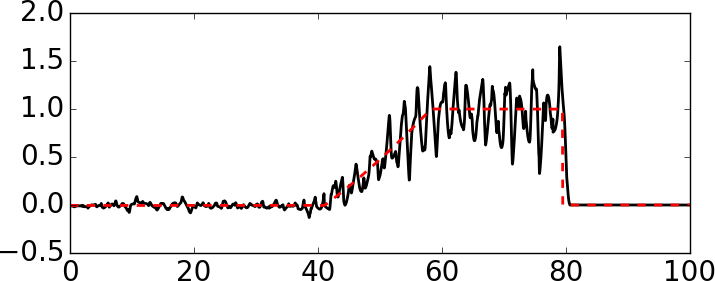}
     \includegraphics[scale=0.23]{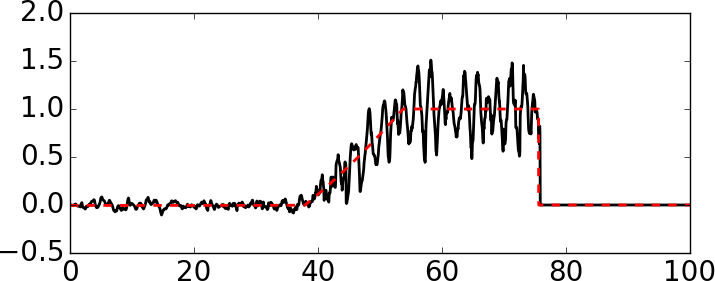}
     \includegraphics[scale=0.23]{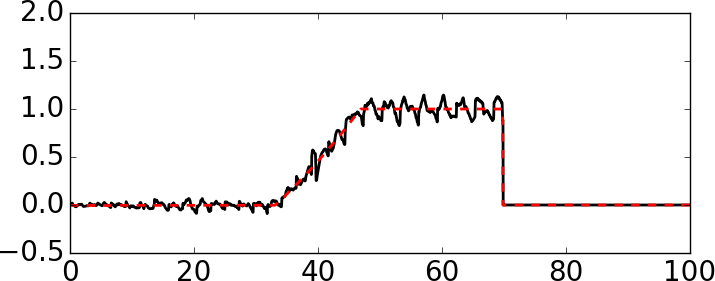}}
   \par\end{centering}
   \caption{\small Effective right-going shock in a layered medium with 
     (a) $\theta=22.5^\circ$, (b) $\theta=45^\circ$ and (c) $\theta=67.5^\circ$.
       The material coefficients are given by \eqref{layered_medium} with
       impedance and sound speed contrast of (from left to right) 
       $Z_B/Z_A=c_B/c_A=4,~2$ and $1.25$ (with $c_A=1$ and $Z_A=1$). 
       The red dashed line is the solution of the homogenized system \eqref{homogenized_leading_order}.
       The initial condition is an effective purely right-going shock
       given by \eqref{init_cond_shock} with
       $x_s=30$, $u_l$ given by \eqref{right-going_shock}, $u_r=0$,
       $\sigma_l=1$ and $\sigma_r=0$.
       We show the solution at $t=30$ and use the nonlinear stress-strain relation \eqref{nonlinearity_exponential}.
     \label{fig:right-going_on_cz-dispersive_media}}
\end{figure}

For $\theta\in(0^\circ, 90^\circ)$, both $C$- and $Z$-dispersion are introduced 
by diffraction and reflections due to variation in the sound speed and impedance, respectively.
In Figure \ref{fig:right-going_on_cz-dispersive_media} we consider 
$\theta=22.5^\circ, ~45^\circ$ and $67.5^\circ$ and media with different sound speed and impedance contrast. 
We use $c_B/c_A=Z_B/Z_A=4, ~2$ and $1.25$ (with $c_A=1$ and $Z_A=1$).
Again, the approximation improves as the sound speed and impedance contrasts are reduced. 

\begin{figure}
 \begin{centering}
   \includegraphics[scale=0.17]{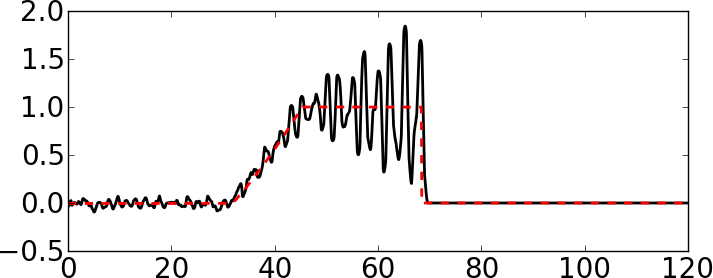}
   \includegraphics[scale=0.17]{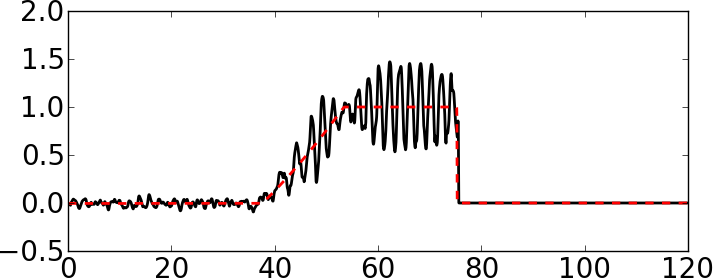}
   \includegraphics[scale=0.17]{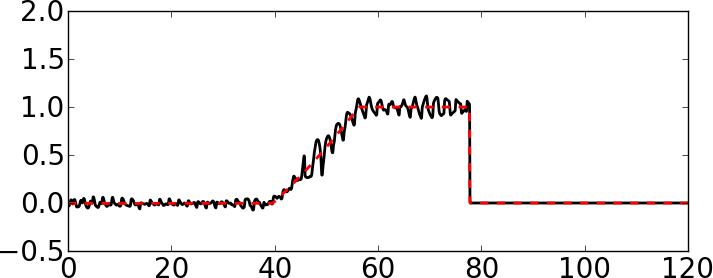}
   \includegraphics[scale=0.17]{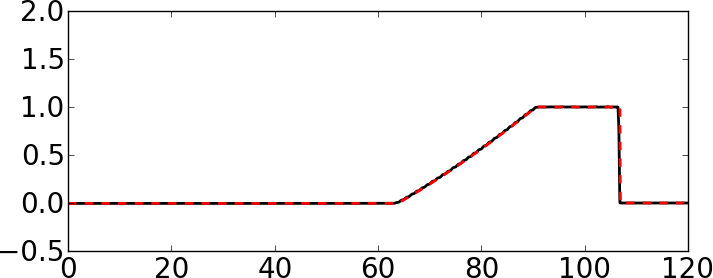}
  \par\end{centering}
  \caption{\small Effective right-going shock in a layered medium with $\theta=90^\circ$.
      The material coefficients are given by \eqref{layered_medium} with
      impedance contrast of (from left to right) $Z_B/Z_A=4,~2,~1.25$ and $1$.
      When $Z_B/Z_A=1$ we use $K_A=\rho_A=1$, $K_B=4$ and $\rho_B=1/4$. 
      For all other cases we use $K_A=\rho_A=1$ and $K_B=\rho_B=Z_B/Z_A$. 
      The red dashed line is the solution of the homogenized system \eqref{homogenized_leading_order}.      
      The initial condition is an effective purely right-going shock
      given by \eqref{init_cond_shock} with
      $x_s=30$, $u_l$ given by \eqref{right-going_shock}, $u_r=0$,
      $\sigma_l=1$ and $\sigma_r=0$.
      We show the solution at $t=40$ and use the nonlinear stress-strain relation \eqref{nonlinearity_exponential}.
  \label{fig:right-going_on_z-dispersive_media}}
\end{figure}

Finally, when $\theta=90^\circ$ $Z$-dispersion is introduced by reflections 
due to variation in the impedance.
We consider media with impedance contrast of $Z_B/Z_A=4,~2,~1.25$ and $1$.
The results are shown in Figure \ref{fig:right-going_on_z-dispersive_media}.
When $Z_B/Z_A=1$ the material parameters are given by 
$K_A=\rho_A=1$, $K_B=4$ and $\rho_B=1/4$; i.e., the material is still non-homogeneous. 
For all other cases we use $K_A=\rho_A=1$ and $K_B=\rho_B=Z_B/Z_A$. 
As before, the leading order approximation improves as the impedance contrast is reduced.

Evidently, the microscopic effects of the induced dispersion (i.e, the oscillatory behavior) cannot
be described by the homogenized system \eqref{homogenized_leading_order}.
This is expected since it is only a leading order approximation. Indeed, one of the assumptions during
the homogenization process is that the wavelength of the pulse is large compared to the
periodicity in the medium,
which is not the case for the right-going shock in these experiments.
Nevertheless,
the homogenized system
seems to accurately predict
not only the speed of the shock but also the macroscopic behavior of the solution.

\begin{remark}[Numerical methods]
  We solve the variable-coefficient system \eqref{p-system} and the homogenized system
  \eqref{homogenized_leading_order} using PyClaw \cite{pyclaw-sisc}.
  Within PyClaw, we use the classic algorithm implemented in Clawpack \cite{clawpack},
  which is a second order method in space and time based 
  on a Lax-Wendroff discretization combined with Total Variation Diminishing (TVD) limiters. 
  This algorithm is a finite volume method for solving hyperbolic conservation laws
  based on the (approximate) solution of Riemann problems. 
  The Riemann solvers can be found in \cite{2014_cylindrical_solitary_waves}.
  The resolution in all the experiments in this work is
  characterized by a mesh size of $\Delta x=\Delta y=1/128$.
\end{remark}

\subsection{Effective shock-speed}
In this section we perform (a total of 1620) numerical simulations to corroborate the
estimate \eqref{speed_of_shock}, which predicts the speed of propagation of shocks in
periodic media.
  The experiments are split into four cases.  For cases
  (a) and (b) we use a layered medium \eqref{layered_medium} with exponential \eqref{nonlinearity_exponential}
  and cubic \eqref{nonlinearity_polynomial} nonlinearities, respectively; for cases
  and (c) and (d) we use a sinusoidal medium \eqref{sin_medium} with exponential \eqref{nonlinearity_exponential}
  and cubic \eqref{nonlinearity_polynomial} nonlinearities, respectively.
For each of these four cases, we consider $\theta\in\{0^\circ, 22.5^\circ, 45^\circ, 67.5^\circ, 90^\circ\}$.
For each value of $\theta$, we perform 81 experiments all starting with a right-going shock
given by \eqref{init_cond_shock}
with $u_l$ given by \eqref{right-going_shock}, $u_r=0$,
$K_A=\rho_A=1$ and all possible combinations of: 
\begin{subequations} \label{param_for_speed_of_shock}
\begin{align}
	\rho_B &= \{2,~3.5,~5\}, &  K_B &= \{2,~3.5,~5\}, \\
	\sigma_l &= \{2,~4,~8\},  & \sigma_r &= \{0,~0.5,~1\}.
\end{align}
\end{subequations}

\begin{figure}
  \begin{centering}
    \includegraphics[scale=0.27]{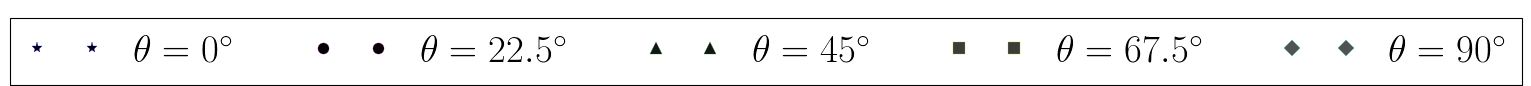}
    \vspace{-10pt}
    
    \subfloat[Layered medium \eqref{layered_medium} with exp. nonlinearity \eqref{nonlinearity_exponential}]
             {\includegraphics[scale=0.25]{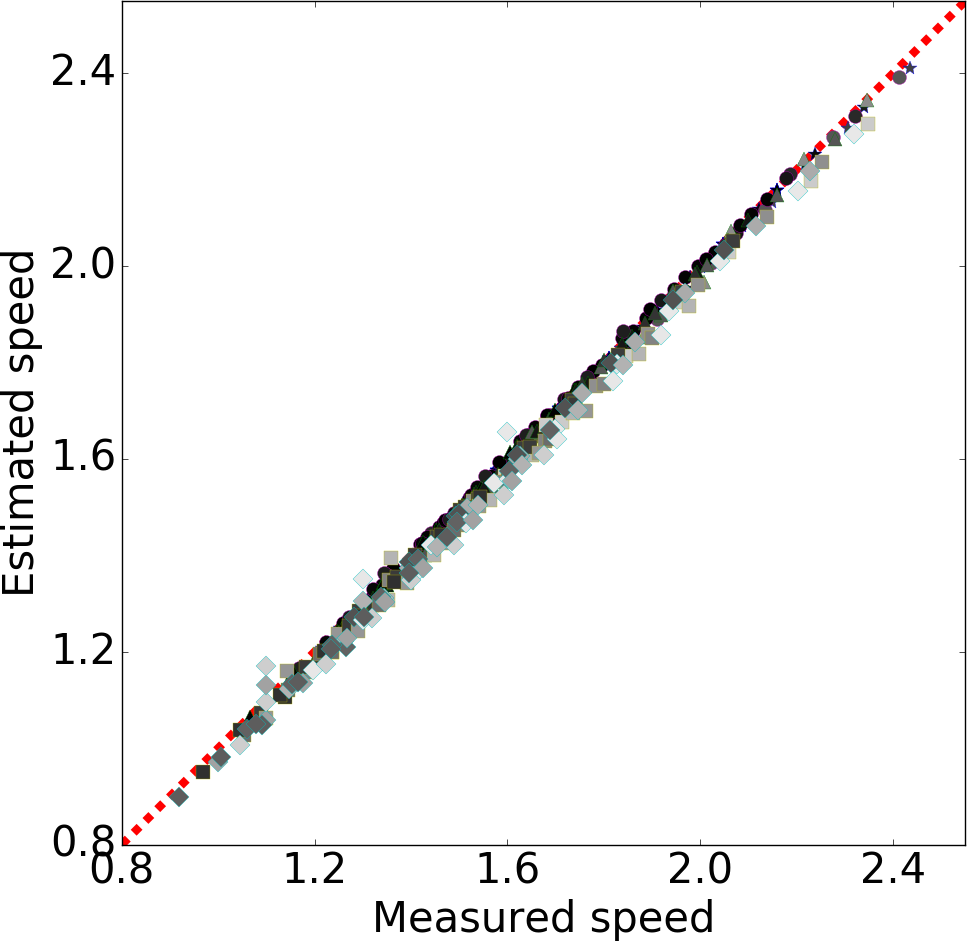}}
    \quad
    \subfloat[Layered medium \eqref{layered_medium} with cubic nonlinearity \eqref{nonlinearity_polynomial}]
             {\includegraphics[scale=0.245]{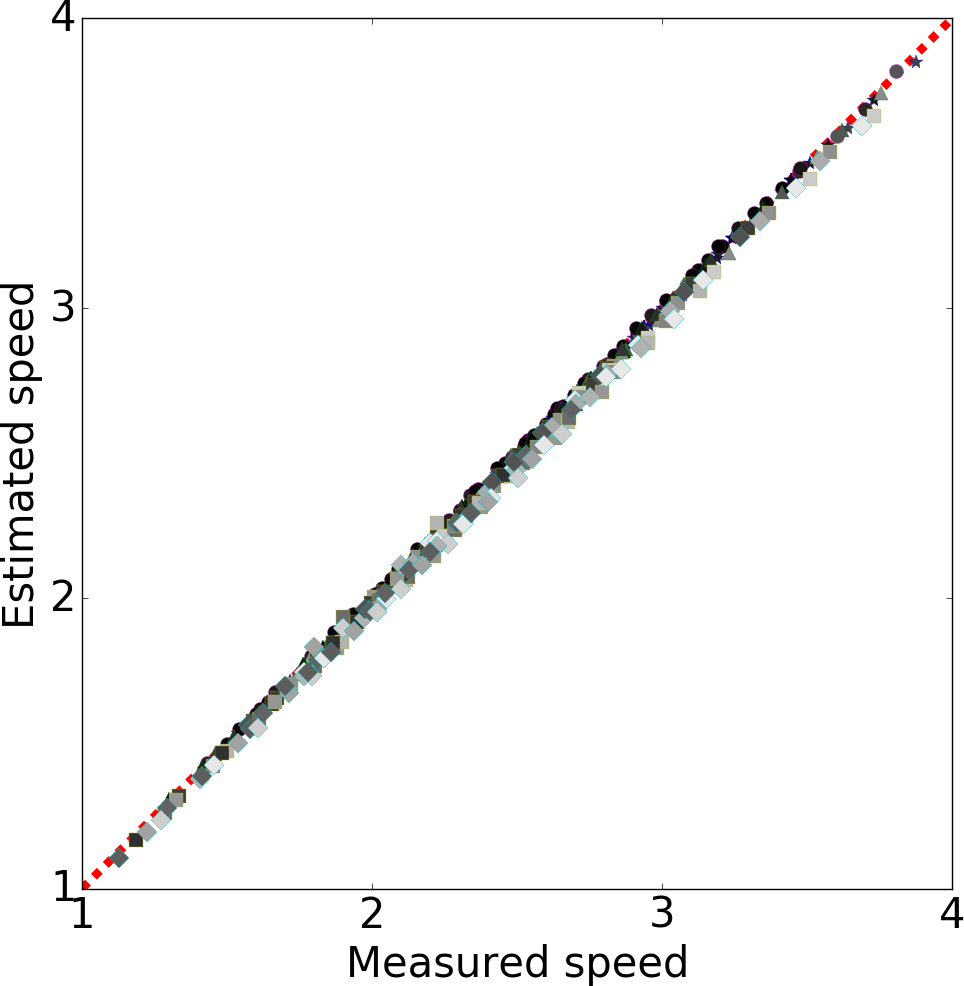}}
             
    \vspace{5pt}
    \subfloat[Sinusoidal medium \eqref{sin_medium} with exp. nonlinearity \eqref{nonlinearity_exponential}]
             {\includegraphics[scale=0.25]{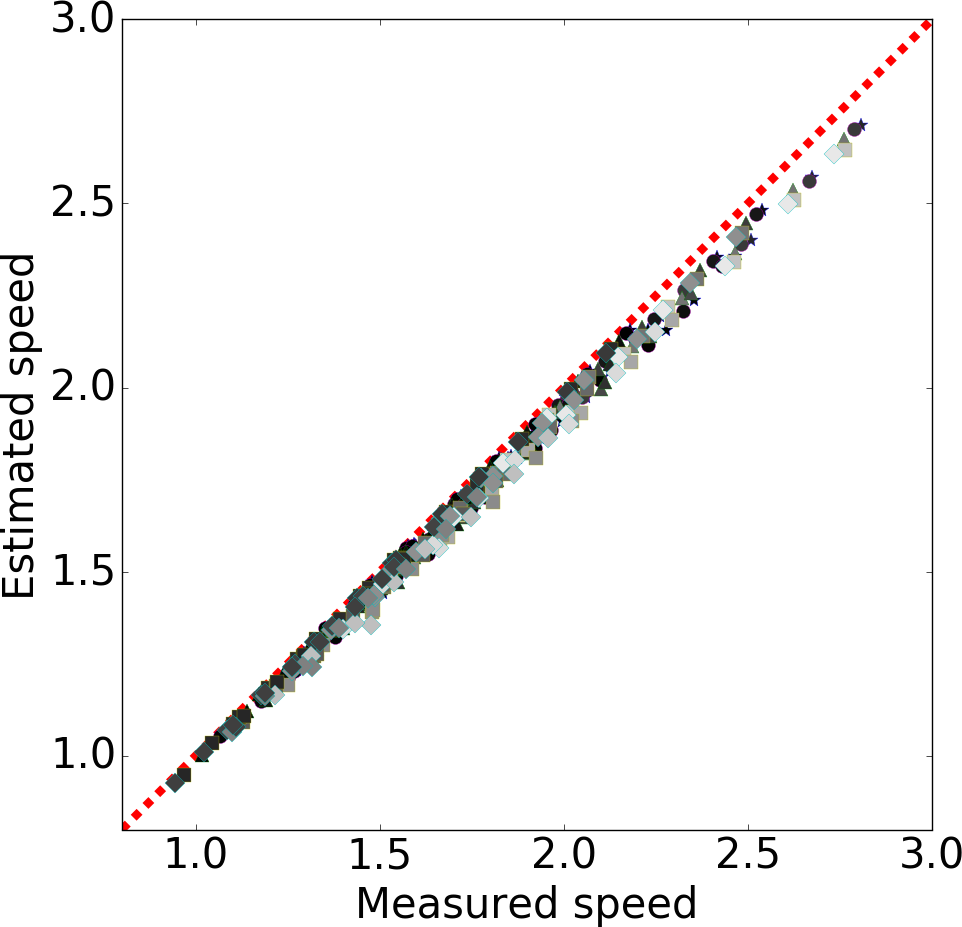}}
    \quad
    \subfloat[Sinusoidal medium \eqref{sin_medium} with cubic nonlinearity \eqref{nonlinearity_polynomial}]
             {\includegraphics[scale=0.26]{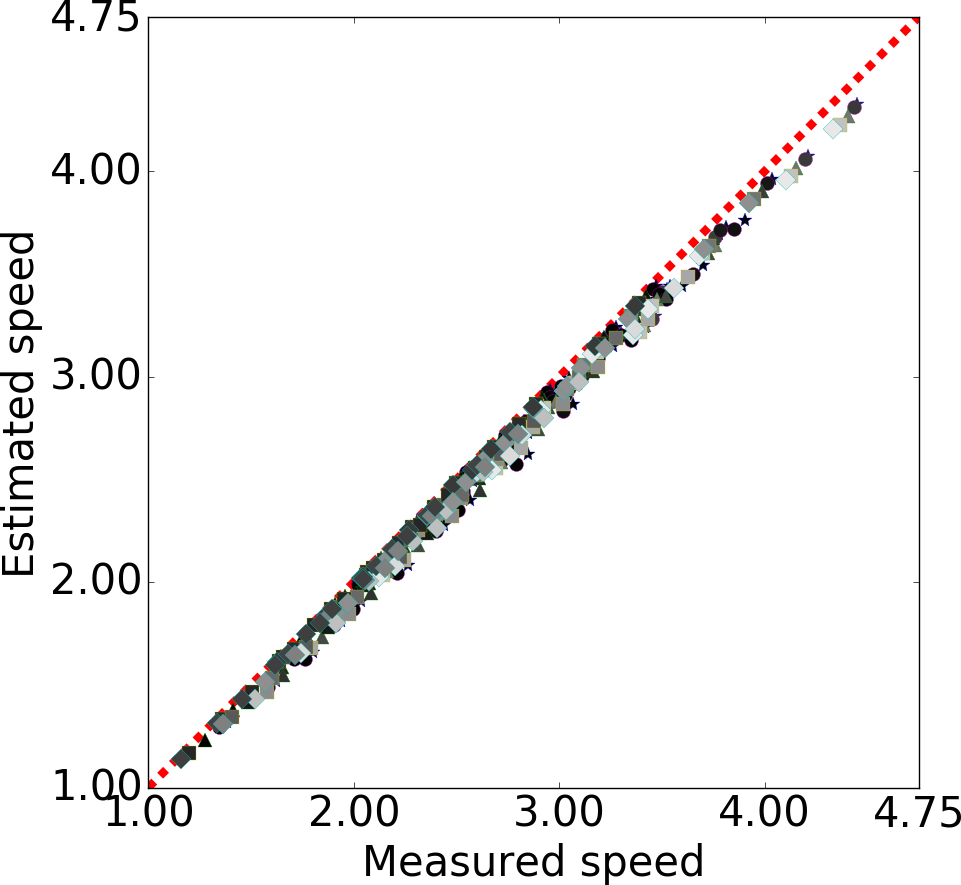}}
    \par\end{centering}
    \caption{\small
      Estimated versus measured speed for 1620 experiments.
      %In the first and second rows the coefficients of the periodic medium are given by the piecewise constant and
      %sinusoidal functions \eqref{layered_medium} and \eqref{sin_medium}, respectively. 
      %
      %The experiments in the left and right columns correspond to the exponential and cubic nonlinearities given by
      %\eqref{nonlinearity_exponential} and \eqref{nonlinearity_polynomial}, respectively. 
      %
      We consider $\theta\in\{0^\circ, 22.5^\circ, 45^\circ, 67.5^\circ, 90^\circ\}$.
      The initial condition is an effective purely right-going shock
      given by \eqref{init_cond_shock} with
      $x_s=40$, $u_l$ given by \eqref{right-going_shock}, $u_r=0$,
      $K_A=\rho_A=1$ and all possible combinations of \eqref{param_for_speed_of_shock}.
      Darker colors in the markers correspond to less dispersive media as
      predicted by \cite[equation (34)]{quezada2014two}.
      \label{fig:estimated_vs_measured}}
\end{figure}

For each experiment, we measure the speed of the shock and compare it with the
speed predicted by \eqref{speed_of_shock}.
Each experiment is represented by a dot in Figure \ref{fig:estimated_vs_measured},
while the predicted speed is indicated by the red dashed line.
The color of the dot represents the amount of dispersion introduced by the medium;
darker colors correspond to less dispersive media
as predicted by \cite[equation (34)]{quezada2014two}.

\section{Towards an effective Lax-entropy condition for periodic media}\label{sec:lax_entropy}
Consider the one-dimensional version of \eqref{p-system} written as a conservation law
and consider only the nonlinear convex stress-strain relation
  \eqref{nonlinearity_exponential}.
A right-going shock for such system satisfies the Lax-entropy condition
\cite{lax1973hyperbolic} if its speed, denoted as $s^+$, satisfies 
\begin{align}
  \lambda^+_l > s^+ > \lambda^+_r,
\end{align}
which imposes that the right-going characteristics (denoted by $\lambda^+$)
from the left and right state of the shock impinge each other.
For a homogeneous medium, this condition is solely dependent upon the solution to the left and right of the shock.

\begin{figure}
  \begin{centering}
    \includegraphics[scale=0.325]{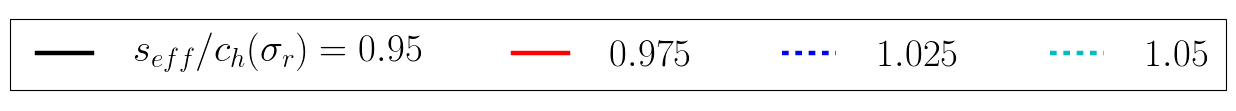}
    \vspace{-15pt}
    
    \subfloat[$C$-dispersive media ($\theta=0^\circ$)\label{fig:entropy_diff_speeds_theta0}]{
      \includegraphics[scale=0.2]{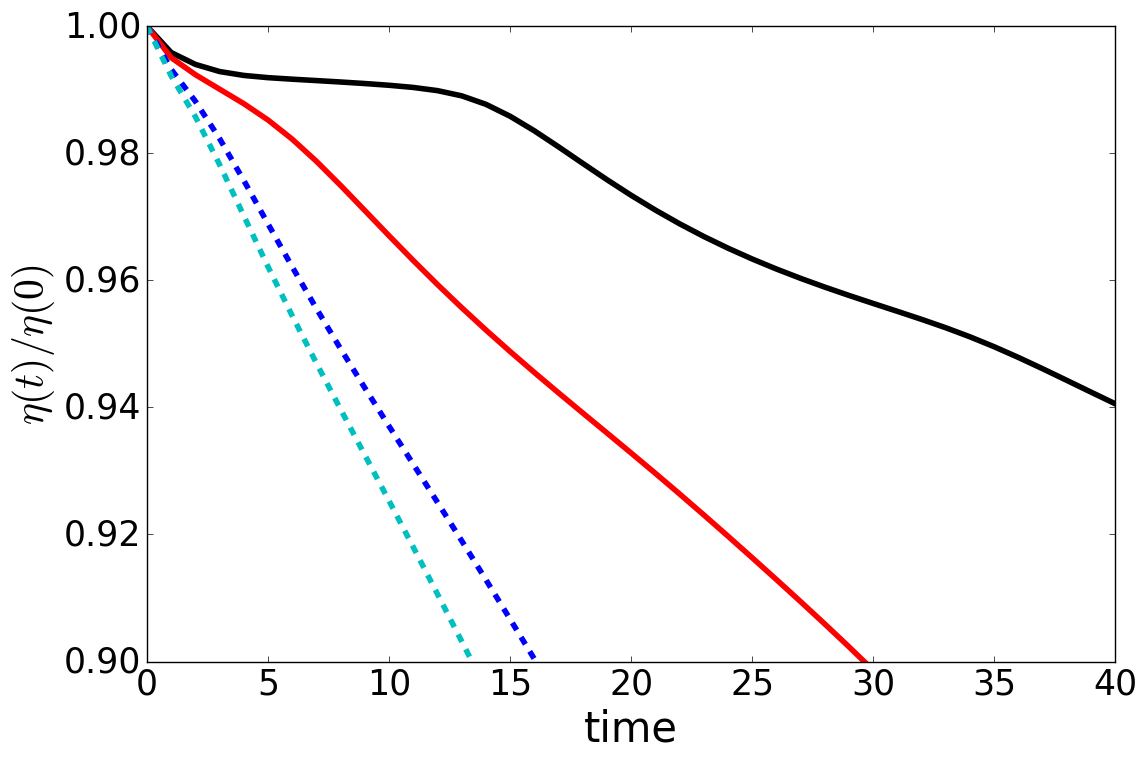}}
    ~
    \subfloat[$Z$-dispersive media ($\theta=90^\circ$)\label{fig:entropy_diff_speeds_theta90}]{
      \includegraphics[scale=0.2]{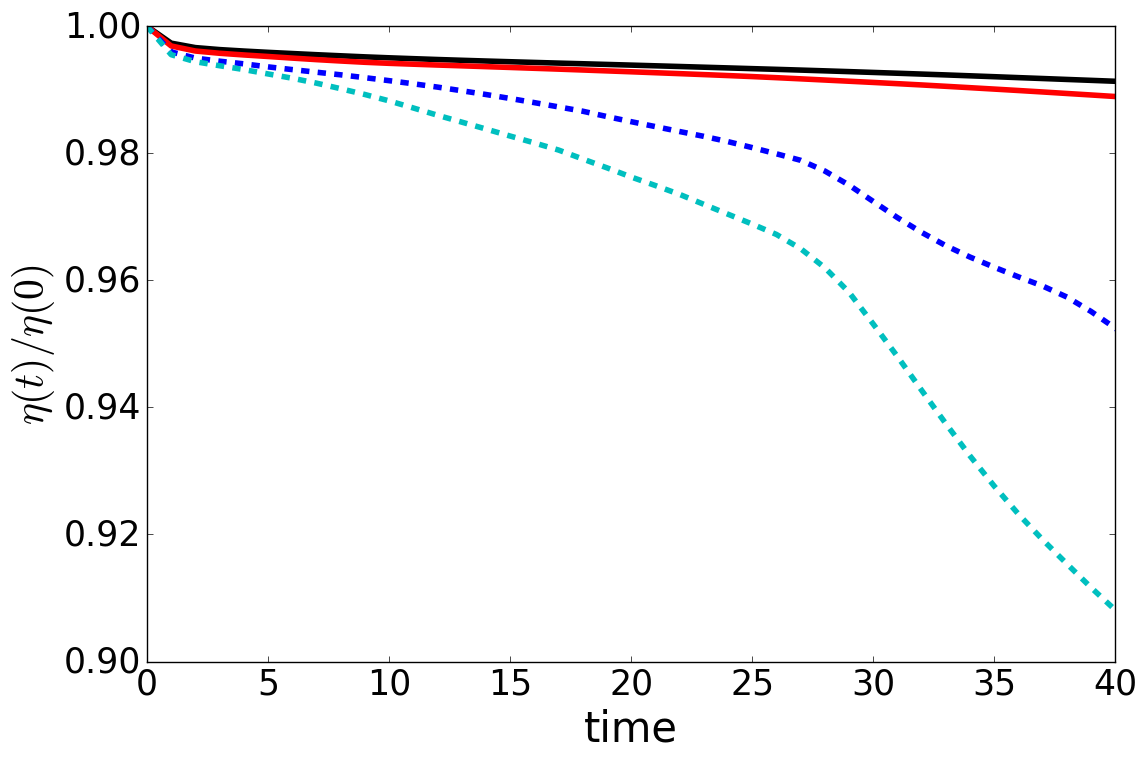}}
    \par\end{centering}
    \caption{\small Normalized entropy evolution for two different media with $\theta=\{0^\circ,~90^\circ\}$.
      The material coefficients are given by the layered medium \eqref{layered_medium} with
      $K_A=\rho_A=\rho_B=1$ and $K_B=16$ $\implies Z_B/Z_A=c_B/c_A=4$.
      The initial condition is an effective purely right-going shock
      given by \eqref{init_cond_shock} with
      $x_s=15$, $u_l$ given by \eqref{right-going_shock}, $u_r=0$,
      $\sigma_r=0$
      and $\sigma_l$ chosen so that (using \eqref{speed_of_shock})
      $s_{\text{eff}}/c_h(\sigma_r)=\{0.95,~0.975,~1.025,~1.05\}$.
      \label{fig:entropy_diff_speeds}}
\end{figure}

In principle, it is reasonable to expect that the effective dispersion introduced in
periodic media
might regularize
weak shocks that can otherwise propagate in homogeneous media.
Therefore, it is desirable to identify a condition for a shock to be able to propagate as a
stable shock in
periodic media. 
This problem is studied in \cite{Ketcheson_LeVeque_2011} for $Z$-dispersive media
(i.e., media as in Figure \ref{fig:layered_intro} with $\theta=90^\circ$).
The authors hypothesize that the speed $s_{\text{eff}}$ of a stable (right-going) shock must satisfy
\begin{align}\label{lax-entropy_z-media}
  s_{\text{eff}} > c_h(\sigma_r) := \left(\int_0^1\left(\frac{\sigma_r^\prime(x)}{\rho(x)}\right)^{-1/2}dx\right)^{-1}.
\end{align}
In other words, a stable right-going shock in a $Z$-dispersive medium must travel faster than the harmonic average of the sound speed.
To corroborate this hypothesis, the authors perform multiple experiments monitoring the evolution of the (global) entropy
\begin{align}\label{entropy}
  \eta(t)= \int_\Omega\left[\frac{1}{2}\rho(x)u^2+\int_0^\epsilon \sigma(z,x)dz\right]dx,
\end{align}
which remains constant for smooth solutions but decreases in the presence of shocks. 
This hypothesis, however, does not seem to accurately predict behavior in the
more general situation depicted in Figure \ref{fig:layered_intro}. 
We illustrate this in Figure \ref{fig:entropy_diff_speeds}, where we plot the evolution of the normalized entropy
$\eta(t)/\eta(0)$ for different experiments.
We consider $\theta\in\{0^\circ,~90^\circ\}$ and $K_A=\rho_A=\rho_B=1$, $K_B=16$.
Note that $Z_B/Z_A=c_B/c_A=4$; i.e, both the impedance and the sound speed change in space.
For each angle $\theta$,
the initial condition is an effective purely right-going shock
given by \eqref{init_cond_shock} with
$x_s=15$, $u_l$ given by \eqref{right-going_shock}, $u_r=0$,
$\sigma_r=0$
and $\sigma_l$ chosen so that (using \eqref{speed_of_shock})
$s_{\text{eff}}/c_h(\sigma_r)=\{0.95,~0.975,~1.025,~1.05\}$.
When $\theta=90^\circ$ (the situation studied in \cite{Ketcheson_LeVeque_2011}), shocks with speed $s_{\text{eff}}>c_h(\sigma_r)$
propagate as stable viscous shocks (i.e., do not get regularized). Otherwise, the shock is regularized by the
induced dispersion and the entropy remains constant. Note that some entropy is lost in all simulations due to numerical dissipation.
The condition \eqref{lax-entropy_z-media} clearly does not hold for $\theta=0^\circ$. 

\begin{figure}
  \begin{centering}
    \includegraphics[scale=0.25]{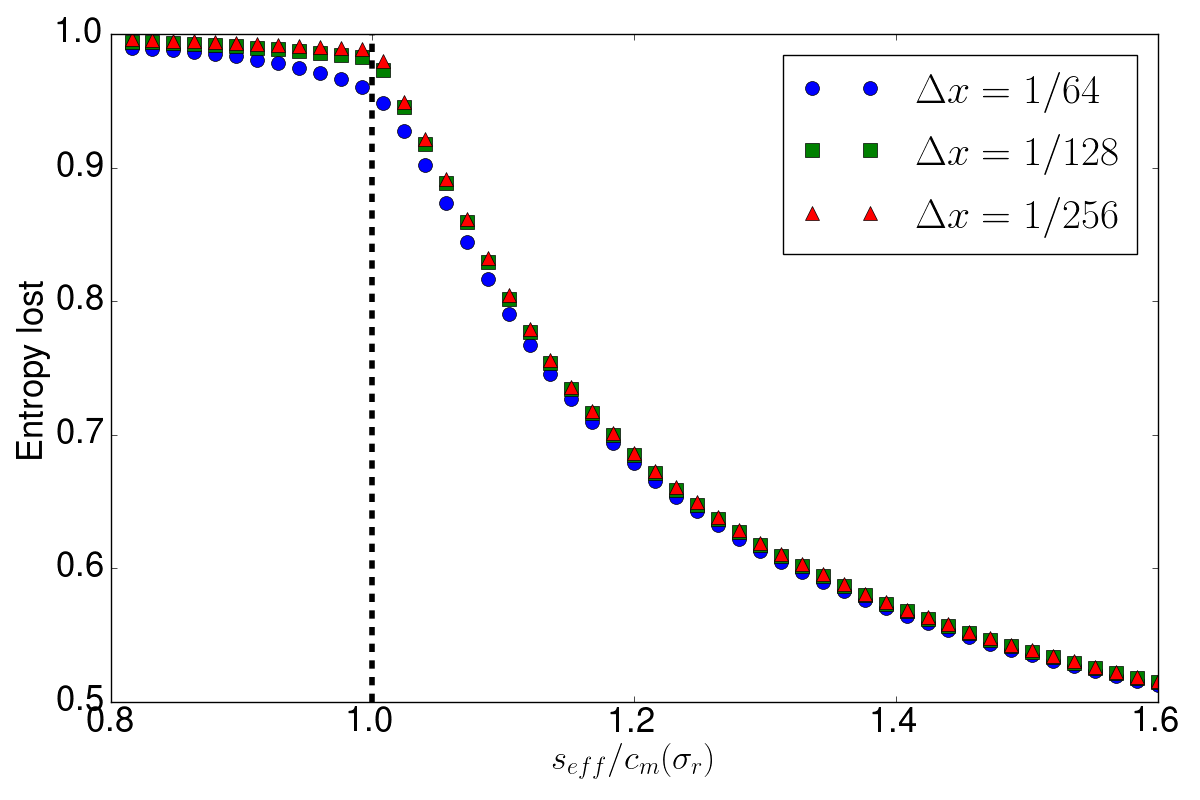}
    \par\end{centering}
    \caption{\small
      Entropy lost at $t=20$ in $C$-dispersive media with $Z=1$.      
      The material coefficients are given by the layered medium \eqref{layered_medium} with
      $K_A=\rho_A=1$, $K_B=4$ and $\rho_B=1/4$ $\implies Z=1, ~c_B/c_A=4$.
      The initial condition is an effective purely right-going shock
      given by \eqref{init_cond_shock} with
      $x_s=15$, $u_l$ given by \eqref{right-going_shock}, $u_r=0$,
      $\sigma_r=0$
      and $\sigma_l$ chosen so that (using \eqref{speed_of_shock})
      $s_{\text{eff}}/c_m(\sigma_r)\in[0.8,~1.6]$.
      \label{fig:entropy_diff_speeds_c-disp}}
\end{figure}

We found empirically that for $C$-dispersive media with $Z=1$, the speed of a stable (right-going) shock satisfies
\begin{align}\label{lax-entropy_c-media}
  s_{\text{eff}}>c_m(\sigma_r):=\int_0^1\frac{\sigma_r^\prime(x)}{\rho(x)}dx.
\end{align}
We illustrate this in Figure \ref{fig:entropy_diff_speeds_c-disp}, where we perform multiple experiments all with
the material parameters given by $K_A=\rho_A=1$, $K_B=4$ and $\rho_B=1/4$ $\implies Z_A=Z_B=1, ~c_B/c_A=4$. 
In all the experiments
the initial condition is an effective purely right-going shock
given by \eqref{init_cond_shock} with
$x_s=15$, $u_l$ given by \eqref{right-going_shock}, $u_r=0$,
$\sigma_r=0$
and $\sigma_l$ chosen so that (using \eqref{speed_of_shock})
$s_{\text{eff}}/c_m(\sigma_r)\in[0.8,~1.6]$.
We measure the entropy lost at $t=20$, which is significant when $s_{\text{eff}}>c_m(\sigma_r)$.
This indicates that viscous shocks remain stable through the entire simulation; otherwise,
the shocks are regularized and the entropy is roughly constant (up to numerical errors).
We consider different mesh widths to demonstrate the effect of numerical dissipation in the entropy lost. 

By performing experiments (not shown here) in more general media (i.e., with $Z\neq 1$) we conclude
that condition \eqref{lax-entropy_c-media} does not hold in general.
%
%Indeed, for $C$-dispersive media we can use the theory in \cite{quezada2014two}
%to infer the dispersion induced by variation in the sound speed. It can be shown that,
%for a fixed sound speed contrast, the dispersive effects vanish as the impedance grows.
%
We are interested in a generalized Lax-entropy condition for
periodic media
(like the one depicted in Figure \ref{fig:layered_intro}); however,
for now that problem remains open.
In \S\ref{sec:conclusions} we discuss potential alternatives to solve this problem. 

\section{Conclusions}\label{sec:conclusions}
Since a shock (or any wave) that propagates in a periodic medium
may be influenced by effective dispersion,
it is important to make the distinction between viscous shocks and dispersive shocks.
Viscous shocks are characterized by the presence of (nearly) discontinuous fronts and by loss in entropy
(which in this work is defined by \eqref{entropy}). 
The main result of this work is an estimate for the speed of a
vanishing-viscosity shock propagating in a
periodic medium,
based on a leading order -- dispersionless -- homogenized system.
Through multiple numerical experiments, we tested the validity of this estimate
and found that the agreement is good. Nevertheless, it is important to remark that the
estimate is expected to hold only when the dispersive effects are relatively small. 
These results represent a generalization of the results presented in \cite{Ketcheson_LeVeque_2011}.

We are also interested in finding a condition for a shock to propagate in
periodic media
without being regularized by the induced
dispersive effects. Such a condition has been developed in \cite{Ketcheson_LeVeque_2011} for the one-dimensional setting.
Although in experiments we have observed similar trend (namely, that larger shocks in
less-dispersive media tend to persist), we have not found an extension of
the specific criterion from \cite{Ketcheson_LeVeque_2011} to the 2D medium considered here.
We found empirically that when $\theta=0^\circ$ and $Z=1$ stable (right-going) shocks
propagate with speed such that $s_{\text{eff}}>c_m(\sigma_r)$ where $c_m(\sigma_r)$ denotes the mean value of the sound speed (to the right of the shock) in the periodic medium.
Nevertheless, this condition does not hold for more general media.
A general condition might be found by deriving a dispersive correction to the
leading order homogenized system \eqref{homogenized_leading_order}
(like the ones obtained in \cite{leveque2003,Fish2001,ketcheson2015diffractons})
and then applying a relaxation method (like in \cite{mazaheri2016first,favrie2017rapid})
that would allow writing the dispersive homogenized system as a hyperbolic model.
In such a model, the characteristic speeds would be modified by the dispersive effects via the relaxation term.
Given such a hyperbolic model one could apply the standard theory of hyperbolic equations to obtain the Lax-entropy condition.
One could also study this problem by deriving the dispersive homogenized system
and then applying the modulation theory of \cite{whitham1965non}.

\section*{Acknowledgment}
This work was supported by funding from King Abdullah University of Science \& Technology (KAUST).

\bibliographystyle{plain}
\bibliography{eff_RH}

\end{document}